\documentclass[12pt]{article}
\usepackage{amsmath}
\usepackage{amsfonts}
\usepackage{amssymb}

\def\thebibliograph#1#2{\section*{{\normalsize \bf #2}}\list
   {[\arabic{enumi}]}{\settowidth\labelwidth{[#1]}\leftmargin\labelwidth
     \advance\leftmargin\labelsep
     \usecounter{enumi}}
     \def\newblock{\hskip .11em plus .33em minus -.07em}
     \sloppy
     \sfcode`\.=1000\relax}

\newtheorem{theorem}{Theorem}
\newtheorem{prop}{Proposition}
\newtheorem{definition}{Definition}

\newtheorem{lemma}{Lemma}
\newtheorem{remark}{Remark}
\input{tcilatex}
\begin{document}

\title{Complex interpolation of variable Triebel-Lizorkin spaces}
\author{Douadi Drihem \\
Department of Mathematics, \\
Laboratory of Functional Analysis and Geometry of spaces, \\
M'sila University, M'sila, Algeria\\
\texttt{douadidr@yahoo.fr}}
\date{\today }
\maketitle

\begin{abstract}
We study complex interpolation of variable Triebel-Lizorkin spaces,
especially we\ present the\ complex interpolation of $F_{p(\cdot
),q}^{\alpha }$ and $F_{p(\cdot ),p(\cdot )}^{\alpha (\cdot )}$ spaces.
Also, some limiting cases are given.

\end{abstract}

\section{Introduction}

Interpolation of spaces have been a central topic in analysis, and are now
of increasing applications in many fields of mathematics especially harmonic
analysis and partial differential equations. For more details on this topic
we refer the reader to Bergh and L\"{o}fstr\"{o}m \cite{BL76}, and Triebel 
\cite{T3}, where the complex interpolation for Besov and/or Triebel-Lizorkin
spaces are given. The main purpose of this paper is to establish the complex
interpolation for variable Triebel-Lizorkin spaces. Firstly we use the
so-called retraction method to present interpolation results in variable
Triebel-Lizorkin spaces $F_{p(\cdot ),q}^{\alpha }$. Secondly\ we shall
apply a method which has been used by \cite{FJ90}\ and\ \cite{SSV13}, where
we shall calculate the Calder\'{o}n products of associated sequence spaces.
Then, from an abstract theory on the relation between the complex
interpolation and the Calder\'{o}n product of Banach lattices obtained by
Calder\'{o}n \cite{Ca64}, Frazier, Jawerth \cite{FJ90}, Mendez, Mitrea \cite%
{MM00} and Kalton, Maybororda, Mitrea \cite{KMM07}, we deduce the complex
interpolation theorems of these sequence spaces. Under some assumptions the
complex interpolation theorems for $F_{p(\cdot ),p(\cdot )}^{\alpha (\cdot
)} $\ are lifted by the $\varphi $-transforms characterization of variable
Triebel-Lizorkin spaces. Finally we will present and briefly discuss some
results concerning the complex interpolation for the spaces $F_{p(\cdot
),q(\cdot )}^{\alpha (\cdot )}$.

As usual, we denote by $\mathbb{R}^{n}$ the $n$-dimensional real Euclidean
space, $\mathbb{N}$ the collection of all natural numbers and $\mathbb{N}%
_{0}:=\mathbb{N}\cup \{0\}$. The letter $\mathbb{Z}$ stands for the set of
all integer numbers.\ The expression $f\lesssim g$ means that $f\leq c\,g$
for some independent constant $c$ (and non-negative functions $f$ and $g$),
and $f\approx g$ means $f\lesssim g\lesssim f$.\vskip5pt

By supp $f$ we denote the support of the function $f$ , i.e., the closure of
its non-zero set. If $E\subset {\mathbb{R}^{n}}$ is a measurable set, then $%
|E|$ stands for the (Lebesgue) measure of $E$ and $\chi _{E}$ denotes its
characteristic function.\vskip5pt

For $v\in \mathbb{Z}$ and $m=(m_{1},...,m_{n})\in \mathbb{Z}^{n}$, let $%
Q_{v,m}$ be the dyadic cube in $\mathbb{R}^{n}$, $Q_{v,m}:=%
\{(x_{1},...,x_{n}):m_{i}\leq 2^{v}x_{i}<m_{i}+1,i=1,2,...,n\}$. For the
collection of all such cubes we use $\mathcal{Q}:=\{Q_{v,m}:v\in \mathbb{Z}%
,m\in \mathbb{Z}^{n}\}$. For each cube $Q$, we denote by $x_{v,m}$ the lower
left-corner $2^{-v}m$ of $Q=Q_{v,m}$ and its side length by $l(Q)$.
Furthermore, we put $v_{Q}:=-\log _{2}l(Q)$, $v_{Q}^{+}:=\max (v_{Q},0)$ and 
$\chi _{Q_{v,m}}=\chi _{v,m},v\in \mathbb{Z},m\in \mathbb{Z}^{n}$.\vskip5pt

The symbol $\mathcal{S}(\mathbb{R}^{n})$ is used in place of the set of all
Schwartz functions on $\mathbb{R}^{n}$.\ We define the Fourier transform of
a function $f\in \mathcal{S}(\mathbb{R}^{n})$ by 
\begin{equation*}
\mathcal{F}(f)(\xi ):=(2\pi )^{-n/2}\int_{\mathbb{R}^{n}}e^{-ix\cdot \xi
}f(x)dx,\quad \xi \in \mathbb{R}^{n}.
\end{equation*}%
We denote by $\mathcal{S}^{\prime }(\mathbb{R}^{n})$ the dual space of all
tempered distributions on $\mathbb{R}^{n}$. For $v\in \mathbb{Z}$, $\varphi
\in \mathcal{S}(\mathbb{R}^{n})$ and $x\in \mathbb{R}^{n}$, we set $%
\widetilde{\varphi }(x):=\overline{\varphi (-x)}$, $\varphi
_{v}(x):=2^{vn}\varphi (2^{v}x)$, and%
\begin{equation*}
\varphi _{v,m}(x):=2^{vn/2}\varphi (2^{v}x-m)=|Q_{v,m}|^{1/2}\varphi
_{v}(x-x_{v,m})\quad \text{if\quad }Q=Q_{v,m}.
\end{equation*}

The variable exponents that we consider are always measurable functions $p$
on $\mathbb{R}^{n}$ with range in $[c,\infty \lbrack $ for some $c>0$. We
denote the set of such functions by $\mathcal{P}_{0}$. The subset of
variable exponents with range $[1,\infty \lbrack $ is denoted by $\mathcal{P}
$. We use the standard notation $p^{-}:$=$\underset{x\in \mathbb{R}^{n}}{%
\text{ess-inf}}$ $p(x)$ and $p^{+}:$=$\underset{x\in \mathbb{R}^{n}}{\text{%
ess-sup }}p(x)$.

The variable exponent modular is defined by $\varrho _{p(\cdot )}(f):=\int_{%
\mathbb{R}^{n}}\rho _{p(x)}(\left\vert f(x)\right\vert )dx$, where $\rho
_{p}(t)=t^{p}$. The variable exponent Lebesgue space $L^{p(\cdot )}$\
consists of measurable functions $f$ on $\mathbb{R}^{n}$ such that $\varrho
_{p(\cdot )}(\lambda f)<\infty $ for some $\lambda >0$. We define the
Luxemburg (quasi)-norm on this space by the formula $\left\Vert f\right\Vert
_{p(\cdot )}:=\inf \big\{\lambda >0:\varrho _{p(\cdot )}\big(\frac{f}{%
\lambda }\big)\leq 1\big\}$. A useful property is that $\left\Vert
f\right\Vert _{p(\cdot )}\leq 1$ if and only if $\varrho _{p(\cdot )}(f)\leq
1$, see \cite{DHHR}, Lemma 3.2.4.

Let $p,q\in \mathcal{P}_{0}$. The Lebesgue-sequence space $L^{p(\cdot
)}(\ell _{q(\cdot )})$ is defined to be the space of all family of functions 
$f_{v},v\geq 0$\ such that%
\begin{equation*}
\big\|\left( f_{v}\right) _{v\geq 0}\big\|_{L^{p(\cdot )}(\ell _{q(\cdot
)})}:=\big\|\big\|\left( f_{v}(x)\right) _{v\geq 0}\big\|_{\ell _{q(x)}}%
\big\|_{p(\cdot )}.
\end{equation*}%
It is easy to show that $L^{p(\cdot )}(\ell _{q(\cdot )})$\ is always a
quasi-normed space\ and it is a normed space, if $\min (p(x),q(x))\geq 1$\
holds point-wise.

We say that $g:\mathbb{R}^{n}\rightarrow \mathbb{R}$ is \textit{locally }log%
\textit{-H\"{o}lder continuous}, abbreviated $g\in C_{\text{loc}}^{\log }$,
if there exists $c_{\log }(g)>0$ such that%
\begin{equation}
\left\vert g(x)-g(y)\right\vert \leq \frac{c_{\log }(g)}{\log
(e+1/\left\vert x-y\right\vert )}  \label{lo-log-Holder}
\end{equation}%
for all $x,y\in \mathbb{R}^{n}$. We say that $g$ satisfies the log\textit{-H%
\"{o}lder decay condition}, if there exists $g_{\infty }\in \mathbb{R}$ and
a constant $c_{\log }>0$ such that%
\begin{equation*}
\left\vert g(x)-g_{\infty }\right\vert \leq \frac{c_{\log }}{\log
(e+\left\vert x\right\vert )}
\end{equation*}%
for all $x\in \mathbb{R}^{n}$. We say that $g$ is \textit{globally}-log%
\textit{-H\"{o}lder continuous}, abbreviated $g\in C^{\log }$, if it is%
\textit{\ }locally log-H\"{o}lder continuous and satisfies the log-H\"{o}%
lder decay\textit{\ }condition.\textit{\ }The constants $c_{\log }(g)$ and $%
c_{\log }$ are called the \textit{locally }log\textit{-H\"{o}lder constant }%
and the log\textit{-H\"{o}lder decay constant}, respectively\textit{.} We
note that all functions $g\in C_{\text{loc}}^{\log }$ always belong to $%
L^{\infty }$. It is known that for $p\in C^{\log }$ we have%
\begin{equation}
\Vert \chi _{B}\Vert _{{p(\cdot )}}\Vert \chi _{B}\Vert _{{p}^{\prime }{%
(\cdot )}}\approx |B|.  \label{DHHR}
\end{equation}%
Also,%
\begin{equation}
\Vert \chi _{B}\Vert _{{p(\cdot )}}\approx |B|^{\frac{1}{p(x)}},\quad x\in B
\label{DHHR1}
\end{equation}%
for small balls $B\subset {\mathbb{R}^{n}}$ ($|B|\leq 2^{n}$), with
constants only depending on the $\log $-H\"{o}lder constant of $p$ (see, for
example, \cite[Section 4.5]{DHHR}). Here ${p}^{\prime }$ denotes the
conjugate exponent of $p$ given by $1/{p(\cdot )}+1/{p}^{\prime }{(\cdot )}%
=1 $. \vskip5pt

Recall that $\eta _{v,m}(x):=2^{nv}(1+2^{v}\left\vert x\right\vert )^{-m}$,
for any $x\in \mathbb{R}^{n}$, $v\in \mathbb{N}_{0}$ and $m>0$. Note that $%
\eta _{v,m}\in L^{1}$ when $m>n$ and that $\big\|\eta _{v,m}\big\|_{1}=c_{m}$
is independent of $v$.

\subsection{Basic tools}

In this subsection we present some results which are useful for us. The
following lemma is from \cite[Lemma 6.1]{DHR}, see also \cite[Lemma 19]%
{KV122}.

\begin{lemma}
\label{DHR-lemma}Let $\alpha \in C_{\mathrm{loc}}^{\log }$ and let $R\geq
c_{\log }(\alpha )$, where $c_{\log }(\alpha )$ is the constant from $%
\mathrm{\eqref{lo-log-Holder}}$ for $\alpha $. Then%
\begin{equation}
2^{v\alpha (x)}\eta _{v,h+R}(x-y)\leq c\text{ }2^{v\alpha (y)}\eta
_{v,h}(x-y)  \label{alpha-est}
\end{equation}%
with $c>0$ independent of $x,y\in \mathbb{R}^{n}$ and $v,h\in \mathbb{N}%
_{0}. $
\end{lemma}

The previous lemma allows us to treat the variable smoothness in many cases
as if it were not variable at all, namely we can move the term inside the
convolution as follows:%
\begin{equation*}
2^{v\alpha (x)}\eta _{v,h+R}\ast \left\vert f\right\vert (x)\leq c\text{ }%
\eta _{v,h}\ast (2^{v\alpha (\cdot )}\left\vert f\right\vert )(x).
\end{equation*}

\begin{lemma}
\label{DHRlemma}Let $p,q\in C^{\log }$ with $1<p^{-}\leq p^{+}<\infty $ and $%
1<q^{-}\leq q^{+}<\infty $. For $m>n$, there exists $c>0$ such that%
\begin{equation*}
\left\Vert (\eta _{v,m}\ast f_{v})_{v}\right\Vert _{L^{p(\cdot )}(\ell
_{q(\cdot )})}\leq c\left\Vert (f_{v})_{v}\right\Vert _{L^{p(\cdot )}(\ell
_{q(\cdot )})}.
\end{equation*}
\end{lemma}

The proof is given in \cite[Theorem 3.2]{DHR}.

Now we introduce the following sequence space.

\begin{definition}
Let $p,q\in \mathcal{P}_{0}$ where $0<p^{+},q^{+}<\infty $ and let $\alpha :%
\mathbb{R}^{n}\rightarrow \mathbb{R}$. Then for all complex valued sequences 
$\lambda :=\{\lambda _{v,m}\}_{v\in \mathbb{N}_{0},m\in \mathbb{Z}%
^{n}}\subset \mathbb{C}$ we define%
\begin{equation*}
f_{p\left( \cdot \right) ,q\left( \cdot \right) }^{\alpha \left( \cdot
\right) }:=\big\{\lambda :\left\Vert \lambda \right\Vert _{f_{p\left( \cdot
\right) ,q\left( \cdot \right) }^{\alpha \left( \cdot \right) }}<\infty %
\big\},
\end{equation*}%
where%
\begin{equation*}
\left\Vert \lambda \right\Vert _{f_{p\left( \cdot \right) ,q\left( \cdot
\right) }^{\alpha \left( \cdot \right) }}:=\Big\|\Big(\sum\limits_{m\in 
\mathbb{Z}^{n}}2^{v(\alpha \left( \cdot \right) +\frac{n}{2})}\lambda
_{v,m}\chi _{v,m}\Big)_{v\geq 0}\Big\|_{L^{p(\cdot )}(\ell _{q(\cdot )})}
\end{equation*}%
and 
\begin{equation*}
f_{\infty ,q}^{\alpha \left( \cdot \right) }:=\big\{\lambda :\left\Vert
\lambda \right\Vert _{f_{\infty ,q}^{\alpha \left( \cdot \right) }}<\infty %
\big\},
\end{equation*}%
where%
\begin{equation*}
\left\Vert \lambda \right\Vert _{f_{\infty ,q}^{\alpha \left( \cdot \right)
}}:=\sup_{Q\in \mathcal{Q}}\frac{1}{|Q|^{1/q}}\Big(\sum%
\limits_{v=v_{Q}^{+}}^{\infty }\int\limits_{Q}\sum\limits_{m\in \mathbb{Z}%
^{n}}2^{v(\alpha \left( x\right) +\frac{n}{2})q}\left\vert \lambda
_{v,m}\right\vert ^{q}\chi _{v,m}(x)dx\Big)^{1/q}.
\end{equation*}
\end{definition}

Notice that the supremum can be taken respect to dyadic cubes with side
length $\leq 1$.

\begin{lemma}
\label{lamda-est}Let $\alpha \in C_{\mathrm{loc}}^{\log }$, $p,q\in C^{\log }
$, $0<p^{+},q^{+}<\infty $, $j\in \mathbb{N}_{0},m\in \mathbb{Z}^{n}$ and $%
x\in Q_{j,m}$. Let $\lambda \in f_{p\left( \cdot \right) ,q\left( \cdot
\right) }^{\alpha \left( \cdot \right) }$. Then there exists $c>0$
independent of $j$ and $m$ such that%
\begin{equation*}
|\lambda _{j,m}|\leq c\text{ }2^{-j(\alpha (x)-\frac{n}{p(x)}+\frac{n}{2}%
)}\left\Vert \lambda \right\Vert _{f_{p\left( \cdot \right) ,q\left( \cdot
\right) }^{\alpha \left( \cdot \right) }}.
\end{equation*}
\end{lemma}

\textbf{Proof.} Let $\lambda \in f_{p\left( \cdot \right) ,q\left( \cdot
\right) }^{\alpha \left( \cdot \right) },j\in \mathbb{N}_{0},m\in \mathbb{Z}%
^{n}$ and $x\in Q_{j,m}$. Using the fact that $2^{j(\alpha \left( x\right)
-\alpha \left( y\right) )}\leq c$ for any $x,y\in Q_{j,m}$, we obtain%
\begin{eqnarray*}
2^{j\alpha \left( x\right) p^{-}}|\lambda _{j,m}|^{p^{-}}
&=&|Q_{j,m}|^{-1}\int_{Q_{j,m}}2^{j\alpha \left( x\right) p^{-}}|\lambda
_{j,m}|^{p^{-}}\chi _{j,m}(y)dy \\
&\leq &c|Q_{j,m}|^{-1}\int_{Q_{j,m}}2^{j\alpha \left( y\right)
p^{-}}|\lambda _{j,m}|^{p^{-}}\chi _{j,m}(y)dy.
\end{eqnarray*}%
Applying H\"{o}lder's inequality to estimate this expression by 
\begin{equation*}
c|Q_{j,m}|^{-1}\big\|2^{j\alpha \left( \cdot \right) p^{-}}|\lambda
_{j,m}|^{p^{-}}\chi _{j,m}\big\|_{p/p^{-}}\left\Vert \chi _{j,m}\right\Vert
_{(p/p^{-})^{\prime }}\leq c\left\Vert \lambda \right\Vert _{f_{p\left(
\cdot \right) ,q\left( \cdot \right) }^{\alpha \left( \cdot \right)
}}^{p^{-}}\left\Vert \chi _{j,m}\right\Vert _{p/p^{-}}^{-1}2^{-j\frac{np^{-}%
}{2}},
\end{equation*}%
where we have used $\mathrm{\eqref{DHHR}}$. Therefore for any $x\in Q_{j,m}$%
\begin{equation*}
|\lambda _{j,m}|\leq c\text{ }2^{-j(\alpha (x)-\frac{n}{p(x)}+\frac{n}{2}%
)}\left\Vert \lambda \right\Vert _{f_{p\left( \cdot \right) ,q\left( \cdot
\right) }^{\alpha \left( \cdot \right) }},
\end{equation*}%
by $\mathrm{\eqref{DHHR1}}$, which completes the proof. \ $\square $

\begin{prop}
\label{prop1}\textit{Let }$\alpha \in C_{\mathrm{loc}}^{\log }$\textit{.
Then }$\lambda =\{\lambda _{v,m}\in \mathbb{C}\}_{v\in \mathbb{N}_{0},m\in 
\mathbb{Z}^{n}}\in f_{\infty ,q}^{\alpha \left( \cdot \right) }$\textit{\ if
and only if for each \ dyadic cube }$Q_{v,m}$ there is a subset $%
E_{Q_{v,m}}\subset Q_{v,m}$ with $|E_{Q_{v,m}}|>|Q_{v,m}|/2$ (or any other,
fixed, number $0<\varepsilon <1$) such that%
\begin{equation*}
\Big\|\Big(\sum_{v=0}^{\infty }\sum\limits_{m\in \mathbb{Z}^{n}}2^{v(\alpha
\left( \cdot \right) +\frac{n}{2})q}|\lambda _{v,m}|^{q}\chi _{E_{v,m}}\Big)%
^{1/q}\Big\|_{\infty }<\infty .
\end{equation*}%
Moreover, the infimum of this expression over all such collections $%
\{E_{Q_{v,m}}\}_{v,m}$ is equivalent to $\left\Vert \lambda \right\Vert
_{f_{\infty ,q}^{\alpha \left( \cdot \right) }}$.
\end{prop}

The proof is given in \cite{D7}. We will make use of the following
statement, see \cite{DHHMS}, Lemma 3.3.

\begin{lemma}
\label{DHHR-estimate}Let $p\in C^{\log }$ with $1\leq p^{-}\leq p^{+}<\infty 
$. Then for every $m>0$ there exists $\gamma =e^{-2m/c_{\log }(1/p)}\in
\left( 0,1\right) $ such that%
\begin{equation*}
\Big(\frac{\gamma }{\left\vert Q\right\vert }\int_{Q}\left\vert
f(y)\right\vert dy\Big)^{p\left( x\right) }\leq \frac{1}{\left\vert
Q\right\vert }\int_{Q}\left\vert f(y)\right\vert ^{p\left( y\right) }dy+\min
\left( \left\vert Q\right\vert ^{m},1\right) g(x)
\end{equation*}%
for every cube $\mathrm{(}$or ball$\mathrm{)}$ $Q\subset \mathbb{R}^{n}$,
all $x\in Q\subset \mathbb{R}^{n}$ and all $f\in L^{p\left( \cdot \right)
}+L^{\infty }$\ with $\left\Vert f\right\Vert _{p\left( \cdot \right)
}+\left\Vert f\right\Vert _{\infty }\leq 1$, where%
\begin{equation*}
g(x):=\left( e+\left\vert x\right\vert \right) ^{-m}+\frac{1}{\left\vert
Q\right\vert }\int_{Q}\left( e+\left\vert y\right\vert \right) ^{-m}dy.
\end{equation*}
\end{lemma}

Notice that in the proof of this theorem we need only that 
\begin{equation*}
\int_{Q}\left\vert f(y)\right\vert ^{p\left( y\right) }dy\leq 1
\end{equation*}%
and/or $\left\Vert f\right\Vert _{\infty }\leq 1$. Moreover if $\left\vert
Q\right\vert \leq 1$, we have $g(x)\leq c$ $\eta _{0,m}(x)$ for any $x\in Q$%
, where $c>0$ is independent of $x$ and $m$.

\section{Variable Triebel-Lizorkin spaces}

The definition of Triebel-Lizorkin spaces of variable smoothness and
integrability is based on the technique of decomposition of unity exactly in
the same manner as in the case of constant exponents. Select a pair of
Schwartz functions $\Phi $ and $\varphi $ satisfy%
\begin{equation}
\text{supp}\mathcal{F}\Phi \subset \overline{B(0,2)}\text{ and }|\mathcal{F}%
\Phi (\xi )|\geq c\text{ if }|\xi |\leq \frac{5}{3}  \label{Ass1}
\end{equation}%
and 
\begin{equation}
\text{supp}\mathcal{F}\varphi \subset \overline{B(0,2)}\backslash B(0,1/2)%
\text{ and }|\mathcal{F}\varphi (\xi )|\geq c\text{ if }\frac{3}{5}\leq |\xi
|\leq \frac{5}{3}  \label{Ass2}
\end{equation}%
where $c>0$. It easy to see that $\int_{\mathbb{R}^{n}}x^{\gamma }\varphi
(x)dx=0$ for all multi-indices $\gamma \in \mathbb{N}_{0}^{n}$.

Now, we define the spaces under consideration.

\begin{definition}
\label{B-F-def}Let $\alpha :\mathbb{R}^{n}\rightarrow \mathbb{R}$ and $%
p,q\in \mathcal{P}_{0}$ with $0<p^{+},q^{+}<\infty $ . Let $\Phi $ and $%
\varphi $ satisfy $\mathrm{\eqref{Ass1}}$ and $\mathrm{\eqref{Ass2}}$,
respectively and we put $\varphi _{v}=2^{vn}\varphi (2^{v}\cdot )$. The
Triebel-Lizorkin space $F_{p(\cdot ),q(\cdot )}^{\alpha (\cdot )}$\ is the
collection of all $f\in \mathcal{S}^{\prime }(\mathbb{R}^{n})$\ such that 
\begin{equation}
\left\Vert f\right\Vert _{F_{p(\cdot ),q(\cdot )}^{\alpha (\cdot )}}:=\big\|%
\left( 2^{v\alpha \left( \cdot \right) }\varphi _{v}\ast f\right) _{v\geq 0}%
\big\|_{L^{p(\cdot )}(\ell _{q(\cdot )})}<\infty ,  \label{B-def}
\end{equation}%
where $\varphi _{0}$ is replaced by $\Phi $.
\end{definition}

The Triebel-Lizorkin spaces with variable smoothness have first been
introduced in \cite{DHR} under much more restrictive conditions on $\alpha
(\cdot )$. If $p(\cdot )$, $q(\cdot )$ and $\alpha (\cdot )$ are constants,
then we derive the well known Triebel-Lizorkin spaces. Taking $\alpha \in 
\mathbb{R}$ and $q\in (0,\infty ]$ as constants we derive the spaces $%
F_{p(\cdot ),q}^{\alpha }$ studied by Xu in \cite{Xu08} and \cite{Xu09}. The
spaces $F_{p(\cdot ),q(\cdot )}^{\alpha (\cdot )}$ are independent of the
particular choice of the system $\left\{ \varphi _{v}\right\} _{v}$\
appearing in their definitions. Moreover, if $\alpha \in C_{\mathrm{loc}%
}^{\log }$ and\ $p,q\in \mathcal{P}_{0}^{\log }$ with $0<p^{-}\leq
p^{+}<\infty $ and $0<q^{-}\leq q^{+}<\infty $, then 
\begin{equation}
\mathcal{S}(\mathbb{R}^{n})\hookrightarrow F_{p(\cdot ),q(\cdot )}^{\alpha
(\cdot )}\hookrightarrow \mathcal{S}^{\prime }(\mathbb{R}^{n}).  \label{emb}
\end{equation}

\begin{definition}
\label{B-F-def2}Let $\alpha :\mathbb{R}^{n}\rightarrow \mathbb{R}$ and $%
0<q<\infty $. Let $\Phi $ and $\varphi $ satisfy $\mathrm{\eqref{Ass1}}$ and 
$\mathrm{\eqref{Ass2}}$, respectively and we put $\varphi _{v}=2^{vn}\varphi
(2^{v}\cdot )$. The Triebel-Lizorkin space $F_{\infty ,q}^{\alpha (\cdot )}$%
\ is the collection of all $f\in \mathcal{S}^{\prime }(\mathbb{R}^{n})$\
such that 
\begin{equation*}
\left\Vert f\right\Vert _{F_{\infty ,q}^{\alpha (\cdot )}}:=\sup_{Q\in 
\mathcal{Q}}\frac{1}{|Q|^{1/q}}\Big(\sum\limits_{v=v_{Q}^{+}}^{\infty
}\int\limits_{Q}2^{v\alpha \left( x\right) q}\left\vert \varphi _{v}\ast
f(x)\right\vert ^{q}dx\Big)^{1/q}<\infty ,
\end{equation*}%
where $\varphi _{0}$ is replaced by $\Phi $.
\end{definition}

For more information about these function spaces, consult \cite{D6} and \cite%
{D7}, with different notation. Notice that the supremum can be taken respect
to dyadic cubes with side length $\leq 1$. One of the key tools to prove the
interpolation property of the spaces is their $\varphi $-transforms
characterization, which transfers the problem from function spaces to their
corresponding sequence spaces. Let $\Phi $ and $\varphi $ satisfy,
respectively $\mathrm{\eqref{Ass1}}$ and $\mathrm{\eqref{Ass2}}$. By \cite[%
pp. 130--131]{FJ90}, there exist \ functions $\Psi \in \mathcal{S}(\mathbb{R}%
^{n})$ satisfying $\mathrm{\eqref{Ass1}}$ and $\psi \in \mathcal{S}(\mathbb{R%
}^{n})$ satisfying $\mathrm{\eqref{Ass2}}$ such that for all $\xi \in 
\mathbb{R}^{n}$%
\begin{equation}
\mathcal{F}\widetilde{\Phi }(\xi )\mathcal{F}\Psi (\xi )+\sum_{j=1}^{\infty }%
\mathcal{F}\widetilde{\varphi }(2^{-j}\xi )\mathcal{F}\psi (2^{-j}\xi
)=1,\quad \xi \in \mathbb{R}^{n}.  \label{Ass4}
\end{equation}

Furthermore, we have the following identity for all $f\in \mathcal{S}%
^{\prime }(\mathbb{R}^{n})$; see \cite[(12.4)]{FJ90}%
\begin{eqnarray*}
f &=&\Psi \ast \widetilde{\Phi }\ast f+\sum_{v=1}^{\infty }\psi _{v}\ast 
\widetilde{\varphi }_{v}\ast f \\
&=&\sum_{m\in \mathbb{Z}^{n}}\widetilde{\Phi }\ast f(m)\Psi (\cdot
-m)+\sum_{v=1}^{\infty }2^{-vn}\sum_{m\in \mathbb{Z}^{n}}\widetilde{\varphi }%
_{v}\ast f(2^{-v}m)\psi _{v}(\cdot -2^{-v}m).
\end{eqnarray*}%
Recall that the $\varphi $-transform $S_{\varphi }$ is defined by setting $%
(S_{\varphi })_{0,m}=\langle f,\Phi _{m}\rangle $ where $\Phi _{m}(x)=\Phi
(x-m)$ and $(S_{\varphi })_{v,m}=\langle f,\varphi _{v,m}\rangle $ where $%
\varphi _{v,m}(x)=2^{vn/2}\varphi (2^{v}x-m)$ and $v\in \mathbb{N}$. The
inverse $\varphi $-transform $T_{\psi }$ is defined by 
\begin{equation*}
T_{\psi }\lambda :=\sum_{m\in \mathbb{Z}^{n}}\lambda _{0,m}\Psi
_{m}+\sum_{v=1}^{\infty }\sum_{m\in \mathbb{Z}^{n}}\lambda _{v,m}\psi _{v,m},
\end{equation*}%
where $\lambda :=\{\lambda _{v,m}\}_{v\in \mathbb{N}_{0},m\in \mathbb{Z}%
^{n}}\subset \mathbb{C}$, see \cite{FJ90}.

Now we present the $\varphi $-transform characterization of these function
spaces, see \cite{DHR} and \cite{YZW151}.

\begin{theorem}
\label{phi-tran}Let $\alpha \in C_{\mathrm{loc}}^{\log }$ and $p,q\in
C^{\log }$ with $0<p^{+},q^{+}<\infty $. \textit{Suppose that }$\Phi $, $%
\Psi \in \mathcal{S}(\mathbb{R}^{n})$ satisfy\ $\mathrm{\eqref{Ass1}}$ and $%
\varphi ,\psi \in \mathcal{S}(\mathbb{R}^{n})$ satisfy\ $\mathrm{\eqref{Ass2}%
}$ such that $\mathrm{\eqref{Ass4}}$ holds. The operators $S_{\varphi
}:F_{p\left( \cdot \right) ,q\left( \cdot \right) }^{\alpha \left( \cdot
\right) }\rightarrow f_{p\left( \cdot \right) ,q\left( \cdot \right)
}^{\alpha \left( \cdot \right) }$ and $T_{\psi }:f_{p\left( \cdot \right)
,q\left( \cdot \right) }^{\alpha \left( \cdot \right) }\rightarrow
F_{p\left( \cdot \right) ,q\left( \cdot \right) }^{\alpha \left( \cdot
\right) }$ are bounded. Furthermore, $T_{\psi }\circ S_{\varphi }$ is the
identity on $F_{p\left( \cdot \right) ,q\left( \cdot \right) }^{\alpha
\left( \cdot \right) }$.
\end{theorem}

Notice that this theorem is true for $F_{\infty ,q}^{\alpha (\cdot )}$
spaces, with $\alpha \in C_{\mathrm{loc}}^{\log }$ and $0<q<\infty $, see 
\cite{D6}.

\section{Complex interpolation}

In this section we study complex interpolation of variable Triebel-Lizorkin
spaces.

\subsection{Complex interpolation for the spaces $F_{p(\cdot ),q}^{\protect%
\alpha }$}

In this subsection we study the complex interpolation of variable
Triebel-Lizorkin spaces $F_{p(\cdot ),q}^{\alpha }$. We use the so-called
retraction method which allows us to reduce the problem to the interpolation
of appropriate sequence spaces, see for instance the monographs [2, 17-18].
More information about complex interpolation of Besov and Triebel-Lizorkin
spaces of fixed exponents can be found in \cite{T3}, \cite{T95}, \cite{YYZ13}
and \cite{WSD15}. See \cite{NS12}\ for the complex interpolation (introduced
by Triebel \cite{T3}) of Besov spaces and Triebel-Lizorkin spaces with
variable exponents. See also \cite{AH} for the complex interpolation of
variable Besov spaces. Complex interpolation between variable Lebesgue
spaces and $BMO$\ (or Hardy spaces) is given in \cite{Ko09}.

Let $A_{0}:=\{z\in \mathbb{C}:0<\func{Re}z<1\}$ and $A:=\{z\in \mathbb{C}%
:0\leq \func{Re}z\leq 1\}$.

\begin{definition}
Let $(X_{0},X_{1})$ be an interpolation couple of Banach lattices. Define $%
\mathcal{F(}X_{0},X_{1})$ as the space of bounded analytic functions $%
g:A_{0}\rightarrow X_{0}+X_{1}$, which extend continuously to the\ closure $A
$, such that the functions $t\rightarrow g(j+it)$ are bounded continuous
functions into $X_{j},j=0,1$, which tend to zero as $\left\vert t\right\vert
\rightarrow \infty $. We endow $\mathcal{F(}X_{0},X_{1})$ with the norm 
\begin{equation*}
\left\Vert g\right\Vert _{\mathcal{F(}X_{0},X_{1})}:=\max \big(%
\sup_{t}\left\Vert g(it)\right\Vert _{X_{0}},\sup_{t}\left\Vert
g(1+it)\right\Vert _{X_{1}}\big).
\end{equation*}%
Further, we define the complex interpolation space 
\begin{equation*}
\lbrack X_{0},X_{1}]_{\theta }:=\{f\in X_{0}+X_{1}:f=g(\theta )\text{ for
some }g\in \mathcal{F(}X_{0},X_{1})\},\quad 0<\theta <1
\end{equation*}%
and 
\begin{equation*}
\left\Vert f\right\Vert _{[X_{0},X_{1}]_{\theta }}:=\inf \{\left\Vert
g\right\Vert _{\mathcal{F(}X_{0},X_{1})}:g\in \mathcal{F(}X_{0},X_{1}),\quad
g(\theta )=f\}.
\end{equation*}
\end{definition}

Let $X$ be a complex Banach space. A function $f:\mathbb{R}^{n}\rightarrow X$
is said to be a simple function if it can be written as%
\begin{equation*}
f=\sum_{j=1}^{N}a_{j}\chi _{A_{j}}
\end{equation*}%
with $a_{j}\in X$ and pairwise disjoint $A_{j}\subset \mathbb{R}^{n}$, $%
\left\vert A_{j}\right\vert <\infty $ $(j=1,...,N)$. Let $p\in \mathcal{P}$.
The Bochner-Lebesgue spaces with variable exponent $L^{p(\cdot )}(\mathbb{R}%
^{n},X)$ is the collection of all measurable functions $f:\mathbb{R}%
^{n}\rightarrow X$ endowed with the norm:%
\begin{equation*}
\left\Vert f\right\Vert _{L^{p(\cdot )}(\mathbb{R}^{n},X)}:=\inf \big\{%
\lambda >0:\varrho _{L^{p(\cdot )}(\mathbb{R}^{n},X)}\big(\frac{f}{\lambda }%
\big)\leq 1\big\}.
\end{equation*}%
The spaces $L^{p(\cdot )}(\mathbb{R}^{n},X)$ have been introduced by C.
Cheng and J. Xu \cite{CX13}.

Let $X_{0}$ and $X_{1}$, be two complex Banach spaces, both linearly and
continuously embedded in a linear complex Hausdorff space $\mathcal{A}$. Two
such Banach spaces are said to be an interpolation couple $(X_{0},X_{1})$.

\begin{theorem}
\label{Lebsgue-int}Let $0<\theta <1$. Let $p_{0},p_{1}\in \mathcal{P}$ with $%
1\leq p_{0}^{+},p_{1}^{+}<\infty $ . We put%
\begin{equation*}
\frac{1}{p(\cdot )}:=\frac{1-\theta }{p_{0}(\cdot )}+\frac{\theta }{%
p_{1}(\cdot )}.
\end{equation*}%
Further let $(X_{0},X_{1})$ be an interpolation couple. Then%
\begin{equation*}
\lbrack L^{p_{0}(\cdot )}(\mathbb{R}^{n},X_{0}),L^{p_{1}(\cdot )}(\mathbb{R}%
^{n},X_{1})]_{\theta }=L^{p(\cdot )}(\mathbb{R}^{n},[X_{0},X_{1}]_{\theta }).
\end{equation*}
\end{theorem}

\textbf{Proof.} Our approach follows essentially \cite{BL76} and \cite{T95}.
It is not hard to see, that the space of simple functions is dense in $%
L^{p_{0}(\cdot )}(\mathbb{R}^{n},X_{0})\cap L^{p_{1}(\cdot )}(\mathbb{R}%
^{n},X_{1})$, and thus also in $L^{p_{0}(\cdot )}(\mathbb{R}%
^{n},X_{0}),L^{p_{1}(\cdot )}(\mathbb{R}^{n},X_{1})]_{\theta }$ by \cite[%
Theorem 4.2.2]{BL76}. From now we consider only simple functions. Let $%
f(x)\neq 0$ be a simple function,%
\begin{equation}
f=\sum_{j=1}^{N}a_{j}\chi _{A_{j}}  \label{simple-fun}
\end{equation}%
with $a_{j}\in X_{0}\cap X_{1}$ and pairwise disjoint $A_{j}\subset \mathbb{R%
}^{n}$, $\left\vert A_{j}\right\vert <\infty $ $(j=1,...,N)$. Let us prove
that%
\begin{equation*}
\left\Vert f\right\Vert _{[L^{p_{0}(\cdot )}(\mathbb{R}^{n},X_{0}),L^{p_{1}(%
\cdot )}(\mathbb{R}^{n},X_{1})]_{\theta }}\leq \left\Vert f\right\Vert
_{L^{p(\cdot )}(\mathbb{R}^{n},[X_{0},X_{1}]_{\theta })}.
\end{equation*}%
By the\ scaling argument, we see that it suffices to consider the case $%
\left\Vert f\right\Vert _{L^{p(\cdot )}(\mathbb{R}^{n},[X_{0},X_{1}]_{\theta
})}=1$. We put for $0\leq \func{Re}z\leq 1$ 
\begin{equation*}
g(x,z)=h(x,z)\left\Vert f(x)\right\Vert _{[X_{0},X_{1}]_{\theta }}^{\big(%
\frac{p(x)}{p_{1}(x)}-\frac{p(x)}{p_{0}(x)}\big)(z-\theta )},
\end{equation*}%
where $h\in \mathcal{F(}X_{0},X_{1})$ for fixed $x\in \mathbb{R}^{n}$, $%
h(x,\theta )=f(x)$ and $h(x,z)=h(y,z)$ if $x,y\in A_{j},j=1,...,N$ and $z$
is defined on the strip $0\leq \func{Re}z\leq 1$, $h(x,z)=0$ for any $x\in 
\mathbb{R}^{n}-\cup _{j=1}^{N}A_{j}$ and 
\begin{equation*}
\left\Vert h(x,it)\right\Vert _{X_{0}},\left\Vert h(x,1+it)\right\Vert
_{X_{1}}\leq \left\Vert f(x)\right\Vert _{[X_{0},X_{1}]_{\theta
}}+\varepsilon ,
\end{equation*}%
$\varepsilon $ is a given positive number. Using this we derive%
\begin{eqnarray*}
\varrho _{L^{p_{0}(\cdot )}(\mathbb{R}^{n},X_{0})}(g(it)) &=&\int_{\mathbb{R}%
^{n}}\left\vert h(x,it)\right\vert ^{p_{0}(x)}\left\Vert f(x)\right\Vert
_{[X_{0},X_{1}]_{\theta }}^{-\big(\frac{p(x)}{p_{1}(x)}-\frac{p(x)}{p_{0}(x)}%
\big)\theta p_{0}(x)}dx \\
&\leq &\int_{\mathbb{R}^{n}}\left\Vert f(x)\right\Vert
_{[X_{0},X_{1}]_{\theta }}^{p(x)}dx+\varepsilon ^{\prime }=1+\varepsilon
^{\prime }.
\end{eqnarray*}%
Similarly,%
\begin{equation*}
\varrho _{L^{p_{1}(\cdot )}(\mathbb{R}^{n},X_{1})}(g(1+it))\leq
1+\varepsilon ^{\prime }.
\end{equation*}%
Thus%
\begin{equation*}
\left\Vert g\right\Vert _{\mathcal{F(}L^{p_{0}(\cdot )}(\mathbb{R}%
^{n},X_{0}),L^{p_{1}(\cdot )}(\mathbb{R}^{n},X_{1}))}\leq 1.
\end{equation*}%
This and $g(\cdot ,\theta )=f(\cdot )$\ we obtain\ 
\begin{equation*}
\left\Vert f\right\Vert _{[L^{p_{0}(\cdot )}(\mathbb{R}^{n},X_{0}),L^{p_{1}(%
\cdot )}(\mathbb{R}^{n},X_{1})]_{\theta }}\leq 1.
\end{equation*}%
Now let us prove that%
\begin{equation*}
\left\Vert f\right\Vert _{L^{p(\cdot )}(\mathbb{R}^{n},[X_{0},X_{1}]_{\theta
})}\leq \left\Vert f\right\Vert _{[L^{p_{0}(\cdot )}(\mathbb{R}%
^{n},X_{0}),L^{p_{1}(\cdot )}(\mathbb{R}^{n},X_{1})]_{\theta }}.
\end{equation*}%
Let $f$ be a simple function of type $\mathrm{\eqref{simple-fun}}$. Further,
let $g(x,z)\in \mathcal{F(}X_{0},X_{1})$ where $g(x,\theta )=f(x)$. We will
prove that 
\begin{equation*}
\left\Vert f\right\Vert _{L^{p(\cdot )}(\mathbb{R}^{n},[X_{0},X_{1}]_{\theta
})}\leq c\left\Vert g\right\Vert _{\mathcal{F(}L^{p_{0}(\cdot )}(\mathbb{R}%
^{n},X_{0}),L^{p_{1}(\cdot )}(\mathbb{R}^{n},X_{1}))}.
\end{equation*}%
From Lemma 4.3.2 in \cite{BL76} and H\"{o}lder's inequality, we obtain%
\begin{eqnarray*}
\left\Vert f(x)\right\Vert _{[X_{0},X_{1}]_{\theta }} &\leq &\Big(\frac{1}{%
1-\theta }\int_{R}\left\Vert g(x,it)\right\Vert _{X_{0}}\mu _{0}(\theta ,t)dt%
\Big)^{1-\theta } \\
&&\times \Big(\frac{1}{\theta }\int_{R}\left\Vert g(x,1+it)\right\Vert
_{X_{1}}\mu _{1}(\theta ,t)dt\Big)^{\theta },
\end{eqnarray*}%
where $\mu _{0}(\theta ,t)$ and $\mu _{1}(\theta ,t)$ are the Poisson
kernels for the strip $A_{0}$, 
\begin{equation*}
\frac{1}{1-\theta }\int_{R}\mu _{0}(\theta ,t)dt=\frac{1}{\theta }%
\int_{R}\mu _{1}(\theta ,t)dt=1.
\end{equation*}%
Using H\"{o}lder's inequality with respect to the variable $t$\ we get 
\begin{eqnarray*}
\left\Vert f(x)\right\Vert _{[X_{0},X_{1}]_{\theta }} &\leq &\Big(\frac{1}{%
1-\theta }\int_{R}\left\Vert g(x,it)\right\Vert _{X_{0}}^{p_{0}(x)}\mu
_{0}(\theta ,t)dt\Big)^{(1-\theta )/p_{0}(x)} \\
&&\times \Big(\frac{1}{\theta }\int_{R}\left\Vert g(x,1+it)\right\Vert
_{X_{1}}^{p_{1}(x)}\mu _{1}(\theta ,t)dt\Big)^{\theta /p_{1}(x)}.
\end{eqnarray*}%
Hence, again by H\"{o}lder's inequality with respect to the variable $x$\ we
get%
\begin{eqnarray*}
\left\Vert f\right\Vert _{L^{p(\cdot )}(\mathbb{R}^{n},[X_{0},X_{1}]_{\theta
})} &\lesssim &\Big\|\Big(\frac{1}{1-\theta }\int_{R}\left\Vert
g(x,it)\right\Vert _{X_{0}}^{p_{0}(\cdot )}\mu _{0}(\theta ,t)dt\Big)%
^{1/p_{0}(\cdot )}\Big\|_{p_{0}(\cdot )}^{1-\theta } \\
&&\times \Big\|\Big(\frac{1}{\theta }\int_{R}\left\Vert g(x,1+it)\right\Vert
_{X_{0}}^{p_{1}(\cdot )}\mu _{0}(\theta ,t)dt\Big)^{1/p_{1}(\cdot )}\Big\|%
_{p_{1}(\cdot )}^{\theta }.
\end{eqnarray*}%
The first norm is bounded by%
\begin{equation*}
\left\Vert g\right\Vert _{\mathcal{F(}L^{p_{0}(\cdot )}(\mathbb{R}%
^{n},X_{0}),L^{p_{1}(\cdot )}(\mathbb{R}^{n},X_{1}))}.
\end{equation*}%
Indeed,%
\begin{eqnarray*}
\int_{\mathbb{R}^{n}}\frac{1}{1-\theta }\int_{R}\Big\|\frac{g(x,it)}{%
\left\Vert g\right\Vert _{\mathcal{F}}}\Big\|_{X_{0}}^{p_{0}(x)}\mu
_{0}(\theta ,t)dtdx &=&\frac{1}{1-\theta }\int_{R}\int_{\mathbb{R}^{n}}\Big\|%
\frac{g(x,it)}{\left\Vert g\right\Vert _{\mathcal{F}}}\Big\|%
_{X_{0}}^{p_{0}(x)}dx\mu _{0}(\theta ,t)dt \\
&\leq &\frac{1}{1-\theta }\int_{R}\mu _{0}(\theta ,t)dt=1.
\end{eqnarray*}%
Similarly for the rest term. The proof of theorem is complete. \ $\square $

Now we present complex interpolation of variable\ Triebel-Lizorkin space $%
F_{p(\cdot ),q}^{\alpha }$. We use the so-called retraction method which
allows us to reduce the problem to the interpolation of appropriate sequence
spaces. We recall that a Banach space $X$ is called a retract of a Banach
space $Y$ if there are linear continuous operators $R:Y\rightarrow X$
(retraction ) and $S:X\rightarrow Y$ (co-retraction ) such that the
composition $RS$ is the identity operator in $X$.

First let us recall the Littlewood-Paley decomposition. Let $\Psi $\ be a
function\ in $\mathcal{S}(\mathbb{R}^{n})$\ satisfying $0\leq \Psi (x)\leq 1$
for all $x\in \mathbb{R}^{n}$, $\Psi (x)=1$\ for\ $\left\vert x\right\vert
\leq 1$\ and\ $\Psi (x)=0$\ for\ $\left\vert x\right\vert \geq 2$.\ We put $%
\mathcal{F}\varphi _{0}(x)=\Psi (x)$, $\mathcal{F}\varphi (x)=\Psi (\frac{x}{%
2})-\Psi (x)$\ and $\mathcal{F}\varphi _{v}(x)=\mathcal{F}\varphi (2^{1-v}x)$
for $v=1,2,3,....$ Then $\{\mathcal{F}\varphi _{v}\}_{v\in \mathbb{N}_{0}}$\
is a resolution of unity, $\sum_{v=0}^{\infty }\mathcal{F}\varphi _{v}(x)=1$
for all $x\in \mathbb{R}^{n}$.\ Thus we obtain the Littlewood-Paley
decomposition 
\begin{equation*}
f=\sum_{v=0}^{\infty }\varphi _{v}\ast f
\end{equation*}%
of all $f\in \mathcal{S}^{\prime }(\mathbb{R}^{n})$ $($convergence in $%
\mathcal{S}^{\prime }(\mathbb{R}^{n}))$.\ We define%
\begin{equation*}
S(f)=(\varphi _{v}\ast f)_{v}
\end{equation*}%
and 
\begin{equation}
R((f_{v})_{v})=\sum_{v=0}^{\infty }\omega _{v}\ast f_{v},  \label{serie}
\end{equation}%
where $\omega _{v}=\varphi _{v-1}+\varphi _{v}+\varphi _{v+1}.$

\begin{theorem}
\label{rect}Let $p\in C^{\log }$ with $1\leq p^{-}<\infty $, $1\leq q\leq
\infty $ and $\alpha \in \mathbb{R}$. Let $\{\mathcal{F}\varphi _{v}\}_{v\in 
\mathbb{N}_{0}}$\ be a resolution of unity. Then $S$ is a co-retraction from 
$F_{p(\cdot ),q}^{\alpha }$into $L^{p(\cdot )}(\ell _{q}^{\alpha })$ and $R$
is a corresponding retraction from $L^{p(\cdot )}(\ell _{q}^{\alpha })$ onto 
$F_{p(\cdot ),q}^{\alpha }$.
\end{theorem}

\textbf{Proof.} The convergence of the series $\mathrm{\eqref{serie}}$ in $%
\mathcal{S}^{\prime }(\mathbb{R}^{n})$ can be obtained by Lemma 4.5 of \cite%
{AH}, where one has to take into consideration $L^{p(\cdot )}(\ell
_{q}^{\alpha })\hookrightarrow \ell _{\infty }^{\alpha }(L^{p(\cdot )})$.
Clearly $R$ is a bounded linear operator from $F_{p(\cdot ),q}^{\alpha }$
into $L^{p(\cdot )}(\ell _{q}^{\alpha })$. Moreover, $RS$ is the identity
operator in $F_{p(\cdot ),q}^{\alpha }$. Using the support properties of $%
\mathcal{F}\varphi _{v}$ and $\mathcal{F}\omega _{v}$, 
\begin{eqnarray*}
\varphi _{v}\ast R((f_{i})_{i}) &=&\sum_{i=v-2}^{v+2}\varphi _{v}\ast \omega
_{i}\ast f_{i} \\
&=&\sum_{k=-2}^{2}\varphi _{v}\ast \omega _{k+v}\ast f_{k+v},\text{ \ \ \ }%
v=0,1,....
\end{eqnarray*}

We can estimate $\left\vert \varphi _{v}\ast \omega _{k+v}\ast
f_{k+v}\right\vert $ by 
\begin{equation*}
(\mathcal{\eta }_{v}\mathcal{\ast }\left\vert f_{k+v}\right\vert ^{t})^{1/t},%
\text{ \ }0<t<\min (p^{-},q),
\end{equation*}%
see \cite{D6}. Applying Lemma \ref{DHRlemma},%
\begin{eqnarray*}
\left\Vert R((f_{i})_{i})\right\Vert _{F_{p(\cdot ),q}^{\alpha }} &\lesssim
&\sum_{k=-2}^{2}\Big\|\left( 2^{v\alpha }f_{k+v}\right) _{v\geq 0}\Big\|%
_{L^{p(\cdot )}(\ell _{q})} \\
&\lesssim &\big\|\left( 2^{v\alpha }f_{v}\right) _{v\geq 0}\big\|%
_{L^{p(\cdot )}(\ell _{q})}.
\end{eqnarray*}%
The proof of theorem is complete. \ $\square $

Now we are ready to formulate the main statement of this subsection.

\begin{theorem}
\label{Inter-rest}Let $0<\theta <1$. Let $p_{0},p_{1}\in C^{\log }$ with $%
1\leq p_{0}^{+},p_{1}^{+}<\infty $, $1\leq q_{0},q_{1}<\infty $ and $\alpha
_{0},\alpha _{1}\in \mathbb{R}$. We put%
\begin{equation*}
\frac{1}{p(\cdot )}:=\frac{1-\theta }{p_{0}(\cdot )}+\frac{\theta }{%
p_{1}(\cdot )},\text{ }\frac{1}{q}:=\frac{1-\theta }{q_{0}}+\frac{\theta }{%
q_{1}}\text{ \ and \ }\alpha :=(1-\theta )\alpha _{0}+\theta \alpha _{1}.
\end{equation*}%
Then%
\begin{equation*}
\lbrack F_{p_{0}(\cdot ),q_{0}}^{\alpha _{0}},F_{p_{1}(\cdot
),q_{1}}^{\alpha _{1}}]_{\theta }=F_{p(\cdot ),q}^{\alpha }
\end{equation*}%
holds in the sense of equivalent norms.
\end{theorem}

\textbf{Proof.} The embedding $\mathrm{\eqref{emb}}$ shows that $%
(F_{p_{0}(\cdot ),q_{0}}^{\alpha _{0}},F_{p_{1}(\cdot ),q_{1}}^{\alpha
_{1}}) $ is an interpolation couple. By Theorem \ref{rect}, 
\begin{equation*}
\left\Vert f\right\Vert _{[F_{p_{0}(\cdot ),q_{0}}^{\alpha
_{0}},F_{p_{1}(\cdot ),q_{1}}^{\alpha _{1}}]_{\theta }}\approx \left\Vert
(\varphi _{v}\ast f)_{v}\right\Vert _{[L^{p_{0}(\cdot )}(\ell
_{q_{0}}^{\alpha _{0}}),L^{p_{1}(\cdot )}(\ell _{q_{1}}^{\alpha
_{1}})]_{\theta }}.
\end{equation*}%
Now by Theorem \ref{Lebsgue-int} the left-hand side is equivalent to $%
\left\Vert (\varphi _{v}\ast f)_{v}\right\Vert _{L^{p(\cdot )}(\ell
_{q}^{\alpha })}$. \ $\square $

\subsection{Complex interpolation for the spaces $F_{p(\cdot ),p(\cdot )}^{%
\protect\alpha (\cdot )}$ and some limiting cases}

In this subsection we present interpolation results in $F_{p(\cdot ),p(\cdot
)}^{\alpha (\cdot )}$\ spaces and some consequences,. We follow the approach
of Frazier and Jawerth \cite{FJ90}, see also \cite{SSV13} and \cite{YYZ13}.
We start by defining the Calder\'{o}n product of two Banach lattices. Let $(%
\mathcal{A},S,%
\mu
)$ be a $\sigma $-finite measure space and let $\mathfrak{M}$ be the class
of all complex-valued, $\mu $-measurable functions on $\mathcal{A}$. Then a
Banach space $X\subset \mathfrak{M}$ is called a Banach lattice of functions
if for every $f\in X$ and $g\in \mathfrak{M}$ with $|g(x)|\leq |f(x)|$ for $%
\mu $-a.e. $x\in X$ one has $g\in X$ and $\left\Vert g\right\Vert _{X}\leq
\left\Vert f\right\Vert _{X}$.

\begin{definition}
Let $(\mathcal{A},S,%
\mu
)$ be a $\sigma $-finite measure space and let $\mathfrak{M}$ be the class
of all complex-valued, $\mu $-measurable functions on $\mathcal{A}$. Suppose
that $X_{0}$ and $X_{1}$ are Banach lattices on $\mathfrak{M}$. Given $%
0<\theta <1$, define the Calder\'{o}n product $X_{0}^{1-\theta }\cdot
X_{1}^{\theta }$ as the collection of all functions $f\in \mathfrak{M}$
satisfying%
\begin{equation*}
\left\Vert f\right\Vert _{X_{0}^{1-\theta }\cdot X_{1}^{\theta }}:=\inf %
\big\{\left\Vert g\right\Vert _{X_{0}}^{1-\theta }\left\Vert h\right\Vert
_{X_{1}}^{\theta }:|f|\leq |g|^{1-\theta }|h|^{\theta },\mu \text{-}a.e.,\
\left\Vert g\right\Vert _{X_{0}}\leq 1,\ \left\Vert h\right\Vert
_{X_{1}}\leq 1\big\}.
\end{equation*}
\end{definition}

\begin{remark}
Calder\'{o}n products have been introduced by Calder\'{o}n \cite{Ca64} (in a
little bit different form which coincides with the above one). Further
properties we refer to, Frazier and Jawerth \cite{FJ90}, Mendez and Mitrea 
\cite{MM00}, Kalton, Mayboroda and Mitrea \cite{KMM07} and Yang, Yuan and
Zhuo \cite{YYZ13}.
\end{remark}

We need a few useful properties, see \cite{YYZ13}.

\begin{lemma}
\label{YYZ-Lemma}Let $(\mathcal{A},S,%
\mu
)$ be a $\sigma $-finite measure space and let $\mathfrak{M}$ be the class
of all complex-valued, $\mu $-measurable functions on $\mathcal{A}$. Suppose
that $X_{0}$ and $X_{1}$ are Banach lattices on $\mathfrak{M}$. Let $%
0<\theta <1$.

$\mathrm{(i)}$ Then the Calder\'{o}n product $X_{0}^{1-\theta }\cdot
X_{1}^{\theta }$ is a Banach space.

$\mathrm{(ii)}$ Define the Calder\'{o}n product $\widetilde{X_{0}^{1-\theta
}\cdot X_{1}^{\theta }}$ as the collection of all functions $f\in \mathfrak{M%
}$ satisfying%
\begin{equation*}
\left\Vert f\right\Vert _{\widetilde{X_{0}^{1-\theta }\cdot X_{1}^{\theta }}%
}:=\inf \big\{M>0:|f|\leq M|g|^{1-\theta }|h|^{\theta },\ \left\Vert
g\right\Vert _{X_{0}}\leq 1,\ \left\Vert h\right\Vert _{X_{1}}\leq 1\big\}.
\end{equation*}%
Then $\widetilde{X_{0}^{1-\theta }\cdot X_{1}^{\theta }}=X_{0}^{1-\theta
}\cdot X_{1}^{\theta }$ follows with equality of norms.
\end{lemma}

Now we turn to the investigation of the Calder\'{o}n products of the
sequence spaces $f_{p(\cdot ),p(\cdot )}^{\alpha (\cdot )}$.

\begin{theorem}
\label{calderon-prod1}Let $0<\theta <1$. Let $p_{0},p_{1}\in C^{\log }$ with 
$1\leq p_{0}^{+},p_{1}^{+}<\infty $ and $\alpha _{0},\alpha _{1}\in C_{%
\mathrm{loc}}^{\log }$. We put%
\begin{equation}
\frac{1}{p(\cdot )}:=\frac{1-\theta }{p_{0}(\cdot )}+\frac{\theta }{%
p_{1}(\cdot )}\text{ and }\alpha (\cdot ):=(1-\theta )\alpha _{0}(\cdot
)+\theta \alpha _{1}(\cdot ).  \label{Th-Cond}
\end{equation}%
Then%
\begin{equation*}
\big(f_{p_{0}(\cdot ),p_{0}(\cdot )}^{\alpha _{0}(\cdot )}\big)^{1-\theta }%
\big(f_{p_{1}(\cdot ),p_{1}(\cdot )}^{\alpha _{1}(\cdot )}\big)^{\theta
}=f_{p(\cdot ),p(\cdot )}^{\alpha (\cdot )}
\end{equation*}%
holds in the sense of equivalent norms.
\end{theorem}

\textbf{Proof}.

\textit{Step 1.} We shall prove%
\begin{equation*}
\big(f_{p_{0}(\cdot ),p_{0}(\cdot )}^{\alpha _{0}(\cdot )}\big)^{1-\theta }%
\big(f_{p_{1}(\cdot ),p_{1}(\cdot )}^{\alpha _{1}(\cdot )}\big)^{\theta
}\hookrightarrow f_{p(\cdot ),p(\cdot )}^{\alpha (\cdot )}.
\end{equation*}%
We suppose, that sequences $\lambda :=(\lambda _{j,m})_{j,m}$, $\lambda
^{i}:=(\lambda _{j,m}^{i})_{j,m},i=0,1$, are given and that%
\begin{equation*}
|\lambda _{j,m}|\leq |\lambda _{j,m}^{0}|^{1-\theta }|\lambda
_{j,m}^{1}|^{\theta }
\end{equation*}%
holds for all $j\in \mathbb{N}_{0}$ and $m\in \mathbb{Z}^{n}$. Let 
\begin{equation}
g(x):=\Big(\sum_{j=0}^{\infty }\sum\limits_{m\in \mathbb{Z}^{n}}2^{j(\alpha
(x)+n/2)p(x)}|\lambda _{j,m}|^{p(x)}\chi _{j,m}(x)\Big)^{1/p(x)}.
\label{notation1}
\end{equation}%
Since, 
\begin{eqnarray*}
&&2^{j(\alpha (x)+n/2)p(x)}|\lambda _{j,m}|^{p(x)}\chi _{j,m}(x) \\
&=&\left( 2^{j(\alpha _{0}(x)+n/2)}|\lambda _{j,m}|\chi _{j,m}(x)\right)
^{p(x)(1-\theta )}\left( 2^{j(\alpha _{1}(x)+n/2)}|\lambda _{j,m}|\chi
_{j,m}(x)\right) ^{p(x)\theta },
\end{eqnarray*}%
then, H\"{o}lder's inequality implies that $g(x)$ can be estimated by%
\begin{eqnarray}
&&\Big(\sum_{j=0}^{\infty }\sum\limits_{m\in \mathbb{Z}^{n}}\left(
2^{j(\alpha _{0}(x)+n/2)}|\lambda _{j,m}|\chi _{j,m}(x)\right) ^{p_{0}(x)}%
\Big)^{(1-\theta )/p_{0}(x)}  \label{Second-est} \\
&&\times \Big(\sum_{j=0}^{\infty }\sum\limits_{m\in \mathbb{Z}^{n}}\left(
2^{j(\alpha _{1}(x)+n/2)}|\lambda _{j,m}|\chi _{j,m}(x)\right) ^{p_{1}(x)}%
\Big)^{\theta /p_{1}(x)}.  \notag
\end{eqnarray}%
Apply H\"{o}lder's inequality again but with conjugate indices $\frac{%
p_{0}(\cdot )}{(1-\theta )}$ and $\frac{p_{1}(\cdot )}{\theta }$, we obtain%
\begin{equation*}
\big\|\lambda \big\|_{f_{p(\cdot ),p(\cdot )}^{\alpha (\cdot )}}\leq \big\|%
\lambda ^{0}\big\|_{f_{p_{0}(\cdot ),p_{0}(\cdot )}^{\alpha _{0}(\cdot
)}}^{1-\theta }\big\|\lambda ^{1}\big\|_{f_{p_{1}(\cdot ),p_{1}(\cdot
)}^{\alpha _{1}(\cdot )}}^{\theta }.
\end{equation*}%
\textit{Step 2.} Now we turn to the proof of%
\begin{equation*}
f_{p(\cdot ),p(\cdot )}^{\alpha (\cdot )}\hookrightarrow \big(f_{p_{0}(\cdot
),p_{0}(\cdot )}^{\alpha _{0}(\cdot )}\big)^{1-\theta }\big(f_{p_{1}(\cdot
),p_{1}(\cdot )}^{\alpha _{1}(\cdot )}\big)^{\theta }.
\end{equation*}%
We will use Lemma \ref{YYZ-Lemma}. Let the sequence $\lambda \in f_{p(\cdot
),p(\cdot )}^{\alpha (\cdot )}$ be given with%
\begin{equation*}
\big\|\lambda \big\|_{f_{p(\cdot ),p(\cdot )}^{\alpha (\cdot )}}\neq 0.
\end{equation*}%
We have to find sequences $\lambda ^{0}$ and $\lambda ^{1}$ such that $%
|\lambda _{j,m}|\leq M|\lambda _{j,m}^{0}|^{1-\theta }|\lambda
_{j,m}^{1}|^{\theta }$ for every $j\in \mathbb{N}_{0}$, $m\in \mathbb{Z}^{n}$
and 
\begin{equation}
\big\|\lambda ^{0}\big\|_{f_{p_{0}(\cdot ),p_{0}(\cdot )}^{\alpha _{0}(\cdot
)}}\lesssim 1\text{ \ and \ }\big\|\lambda ^{1}\big\|_{f_{p_{1}(\cdot
),p_{1}(\cdot )}^{\alpha _{1}(\cdot )}}\lesssim 1  \label{Est-Int}
\end{equation}%
with some constant $c$ independent of $\lambda $. We follow ideas of the
proof of Theorem 8.2 in Frazier and Jawerth \cite{FJ90}, see also Sickel,
Skrzypczak and Vyb\'{\i}ral \cite{SSV13}. Set%
\begin{equation*}
u(\cdot ):=p(\cdot )\theta \big(\tfrac{\alpha _{1}(\cdot )}{p_{0}(\cdot )}-%
\tfrac{\alpha _{0}(\cdot )}{p_{1}(\cdot )}\big)+\tfrac{n}{2}\big(\tfrac{%
p(\cdot )}{p_{0}(\cdot )}-1\big)
\end{equation*}%
and 
\begin{equation*}
v(\cdot ):=p(\cdot )(1-\theta )\big(\tfrac{\alpha _{0}(\cdot )}{p_{1}(\cdot )%
}-\tfrac{\alpha _{1}(\cdot )}{p_{0}(\cdot )}\big)+\tfrac{n}{2}\big(\tfrac{%
p(\cdot )}{p_{1}(\cdot )}-1\big).
\end{equation*}%
We put%
\begin{equation*}
\lambda _{j,m}^{0}:=2^{ju(x_{j,m})}\Big(\tfrac{|\lambda _{j,m}|}{\big\|%
\lambda \big\|_{f_{p(\cdot ),p(\cdot )}^{\alpha (\cdot )}}}\Big)%
^{p(x_{j,m})/p_{0}(x_{j,m})},\quad x_{j,m}=2^{-j}m.
\end{equation*}%
Also, set%
\begin{equation*}
\lambda _{j,m}^{1}:=2^{jv(x_{j,m})}\Big(\tfrac{|\lambda _{j,m}|}{\big\|%
\lambda \big\|_{f_{p(\cdot ),p(\cdot )}^{\alpha (\cdot )}}}\Big)%
^{p(x_{j,m})/p_{1}(x_{j,m})}.
\end{equation*}%
Observe that%
\begin{equation*}
|\lambda _{j,m}|\leq \big\|\lambda \big\|_{f_{p(\cdot ),p(\cdot )}^{\alpha
(\cdot )}}\big(\lambda _{j,m}^{0}\big)^{1-\theta }\big(\lambda _{j,m}^{1}%
\big)^{\theta },
\end{equation*}%
which holds now for all pairs $(j,m)$.

\textbf{Estimation of }$\big\|\lambda ^{0}\big\|_{f_{p_{0}(\cdot
),p_{0}(\cdot )}^{\alpha _{0}(\cdot )}}$. Set%
\begin{equation*}
I(\cdot ):=\sum_{j=0}^{\infty }\sum\limits_{m\in \mathbb{Z}^{n}}\Big(%
2^{j(\alpha _{0}(\cdot )+\frac{n}{2})}|\lambda _{j,m}^{0}|\chi _{j,m}\Big)%
^{p_{0}(\cdot )}.
\end{equation*}%
We use the local log-H\"{o}lder continuity of $\alpha _{0},\alpha _{1},p_{0}$
and $p_{1}$ to show that%
\begin{equation}
2^{ju(x)}\leq c\text{ }2^{ju(y)}\quad \text{and}\quad 2^{j(\alpha (x)-\frac{n%
}{p(x)})}\leq c\text{ }2^{j(\alpha (y)-\frac{n}{p(y)})}\text{, \ \ }x,y\in
Q_{j,m},  \label{est-lamda}
\end{equation}%
where $c>0$ is independent of $m\in \mathbb{Z}^{n}$ and $j\in \mathbb{N}_{0}$%
. Now $\Big(\tfrac{|\lambda _{j,m}|}{\big\|\lambda \big\|_{f_{p(\cdot
),p(\cdot )}^{\alpha (\cdot )}}}\Big)^{\frac{p(x_{j,m})}{p_{0}(x_{j,m})}}$\
can be rewritten us%
\begin{equation*}
\Big(2^{j(\alpha (x_{j,m})+\frac{n}{2}-\frac{n}{p(x_{j,m})})}\tfrac{|\lambda
_{j,m}|}{\big\|\lambda \big\|_{f_{p(\cdot ),p(\cdot )}^{\alpha (\cdot )}}}%
\Big)^{\frac{p(x_{j,m})}{p_{0}(x_{j,m})}}2^{-j(\alpha (x_{j,m})+\frac{n}{2}-%
\frac{n}{p(x_{j,m})})\frac{p(x_{j,m})}{p_{0}(x_{j,m})}}.
\end{equation*}%
After applying Lemma \ref{lamda-est}, we get 
\begin{equation}
2^{j(\alpha (x_{j,m})+\frac{n}{2}-\frac{n}{p(x_{j,m})})}\tfrac{|\lambda
_{j,m}|}{\big\|\lambda \big\|_{f_{p(\cdot ),p(\cdot )}^{\alpha (\cdot )}}}%
\lesssim 1.  \label{lemma-cond}
\end{equation}%
Observe that%
\begin{eqnarray*}
2^{j(\alpha (x_{j,m})+\frac{n}{2}-\frac{n}{p(x_{j,m})})}\tfrac{|\lambda
_{j,m}|}{\big\|\lambda \big\|_{f_{p(\cdot ),p(\cdot )}^{\alpha (\cdot )}}}
&=&\frac{1}{\left\vert Q_{j,m}\right\vert }\int_{Q_{j,m}}2^{j(\alpha
(x_{j,m})+\frac{n}{2}-\frac{n}{p(x_{j,m})})}\tfrac{|\lambda _{j,m}|}{\big\|%
\lambda \big\|_{f_{p(\cdot ),p(\cdot )}^{\alpha (\cdot )}}}dy \\
&\lesssim &\frac{1}{\left\vert Q_{j,m}\right\vert }\int_{Q_{j,m}}2^{j(\alpha
(y)+\frac{n}{2}-\frac{n}{p(y)})}\tfrac{|\lambda _{j,m}|}{\big\|\lambda \big\|%
_{f_{p(\cdot ),p(\cdot )}^{\alpha (\cdot )}}}dy
\end{eqnarray*}%
for any $m\in \mathbb{Z}^{n}$ and any $j\in \mathbb{N}_{0}$. Taking the $%
\frac{\sigma (x_{j,m})p(x_{j,m})}{p_{0}(x_{j,m})}$-power, with $\frac{%
p_{0}(\cdot )}{p(\cdot )}<\sigma (\cdot )<1$ and $\sigma \in C^{\log }$, by $%
\mathrm{\eqref{lemma-cond}}$ we can apply Lemma \ref{DHHR-estimate} and
obtain that%
\begin{eqnarray*}
&&\Big(2^{j(\alpha (x_{j,m})+\frac{n}{2}-\frac{n}{p(x_{j,m})})}\tfrac{%
|\lambda _{j,m}|}{\big\|\lambda \big\|_{f_{p(\cdot ),p(\cdot )}^{\alpha
(\cdot )}}}\Big)^{\frac{\sigma (x_{j,m})p(x_{j,m})}{p_{0}(x_{j,m})}} \\
&\lesssim &\frac{1}{\left\vert Q_{j,m}\right\vert }\int_{Q_{j,m}}\Big(%
2^{j(\alpha (y)+\frac{n}{2}-\frac{n}{p(y)})}\tfrac{|\lambda _{j,m}|}{\big\|%
\lambda \big\|_{f_{p(\cdot ),p(\cdot )}^{\alpha (\cdot )}}}\Big)^{\frac{%
\sigma (y)p(y)}{p_{0}(y)}}dy+2^{-jnh}\eta _{0,h}(x_{j,m}).
\end{eqnarray*}%
Applying again Lemma \ref{lamda-est}, the last term with $\frac{1}{\sigma
(x_{j,m})}$-power is bounded by%
\begin{eqnarray*}
&&\frac{2^{\frac{1}{\sigma ^{-}}-1}}{\left\vert Q_{j,m}\right\vert }%
\int_{Q_{j,m}}\Big(2^{j(\alpha (y)+\frac{n}{2}-\frac{n}{p(y)})}\tfrac{%
|\lambda _{j,m}|}{\big\|\lambda \big\|_{f_{p(\cdot ),p(\cdot )}^{\alpha
(\cdot )}}}\Big)^{\frac{p(y)}{p_{0}(y)}}dy+2^{\frac{1}{\sigma ^{-}}%
-jnh+1}\eta _{0,h}(x_{j,m}) \\
&\lesssim &\eta _{j,\rho }\ast \Big(2^{j(\alpha (\cdot )+\frac{n}{2}-\frac{n%
}{p(\cdot )})}\tfrac{|\lambda _{j,m}|}{\big\|\lambda \big\|_{f_{p(\cdot
),p(\cdot )}^{\alpha (\cdot )}}}\chi _{Q_{j,m}}\Big)^{\frac{p(\cdot )}{%
p_{0}(\cdot )}}(x)+2^{-jnh}\eta _{0,h}(x)
\end{eqnarray*}%
for any $h,\rho >0$ and any $x\in Q_{j,m}$. Applying the second estimate of $%
\mathrm{\eqref{est-lamda}}$, the term $\Big(\tfrac{|\lambda _{j,m}|}{\big\|%
\lambda \big\|_{f_{p(\cdot ),p(\cdot )}^{\alpha (\cdot )}}}\Big)^{\frac{%
p(x_{j,m})}{p_{0}(x_{j,m})}}$ can be estimated by%
\begin{equation*}
c\text{ }\eta _{j,\rho _{1}}\ast \Big(\tfrac{|\lambda _{j,m}|}{\big\|\lambda %
\big\|_{f_{p(\cdot ),p(\cdot )}^{\alpha (\cdot )}}}\chi _{Q_{j,m}}\Big)^{%
\frac{p(\cdot )}{p_{0}(\cdot )}}(x)+c\text{ }2^{-j(nh+d)}\eta _{0,h}(x),
\end{equation*}%
where%
\begin{equation*}
d=\big((\alpha +\frac{n}{2}-\frac{n}{p})\frac{p}{p_{0}}\big)^{-}\text{ and }%
\rho _{1}=\rho -c_{\log }((\alpha +\frac{n}{2}-\frac{n}{p})\frac{p}{p_{0}}).
\end{equation*}%
From the estimations above, for any $x\in Q_{j,m}$%
\begin{equation*}
2^{j(\alpha _{0}(x)+\frac{n}{2})}|\lambda _{j,m}^{0}|\lesssim \eta _{j,\rho
_{2}}\ast \digamma _{j,m}(x)+2^{-j(nh+d-(u^{+}+\alpha _{0}^{+}+\frac{n}{2}%
))}\eta _{0,h}(x),
\end{equation*}%
where the implicit positive constant not depending on $x,m$ and $j$, with 
\begin{equation*}
\digamma _{j,m}=2^{j(\alpha _{0}(\cdot )+\frac{n}{2})+ju(\cdot )}\Big(\tfrac{%
|\lambda _{j,m}|}{\big\|\lambda \big\|_{f_{p(\cdot ),p(\cdot )}^{\alpha
(\cdot )}}}\Big)^{\frac{p(\cdot )}{p_{0}(\cdot )}}\chi _{Q_{j,m}},
\end{equation*}%
where 
\begin{equation*}
\rho _{2}=\rho _{1}-c_{\log }(\alpha _{0})-c_{\log }(u).
\end{equation*}%
Hence $\big\|\left( I^{p_{0}(\cdot )}(\cdot )\right) ^{1/p_{0}(\cdot )}\big\|%
_{p_{0}(\cdot )}$ is bounded by%
\begin{equation*}
c\Big\|\Big(\sum_{j=0}^{\infty }\sum\limits_{m\in \mathbb{Z}^{n}}\left( \eta
_{j,\rho _{2}}\ast \digamma _{j,m}\right) ^{p_{0}(\cdot )}\Big)%
^{1/p_{0}(\cdot )}\Big\|_{p_{0}(\cdot )}+c\big\|\eta _{0,h}\big\|%
_{p_{0}(\cdot )}\Big(\sum_{j=0}^{\infty }2^{-j(nh+d-(u^{+}+\alpha _{0}^{+}+%
\frac{n}{2}))p_{0}^{-}}\Big)^{\frac{1}{p_{0}^{-}}},
\end{equation*}%
where in the second estimate we used the embedding $L^{p(\cdot )}(\ell
^{p^{-}})\hookrightarrow L^{p(\cdot )}(\ell ^{p(\cdot )})$. By taking $h$
large enough such that $h>(\alpha _{0}^{+}+u^{+}+\frac{n}{2})/n-\frac{d}{n}$
the second term is bounded. Taking $\rho _{2}$ large enough and applying
Lemma\ \ref{DHRlemma} to estimate the first expression by%
\begin{eqnarray*}
&&c\Big\|\Big(\sum_{j=0}^{\infty }\sum\limits_{m\in \mathbb{Z}^{n}}\digamma
_{j,m}^{p_{0}(\cdot )}\Big)^{1/p_{0}(\cdot )}\Big\|_{p_{0}(\cdot )} \\
&=&c\Big\|\Big(\sum_{j=0}^{\infty }\sum\limits_{m\in \mathbb{Z}%
^{n}}2^{j(\alpha (\cdot )+\frac{n}{2})p(\cdot )}\Big(\tfrac{|\lambda _{j,m}|%
}{\big\|\lambda \big\|_{f_{p(\cdot ),p(\cdot )}^{\alpha (\cdot )}}}\Big)%
^{q(\cdot )}\chi _{Q_{j,m}}\Big)^{1/p_{0}(\cdot )}\Big\|_{p_{0}(\cdot )},
\end{eqnarray*}%
since $u(\cdot )+\alpha _{0}(\cdot )=\alpha (\cdot )\frac{p(\cdot )}{%
p_{0}(\cdot )}+\frac{n}{2}\big(\frac{p(\cdot )}{p_{0}(\cdot )}-1\big)$.
Therefore, $\big\|\left( I^{p_{0}(\cdot )}(\cdot )\right) ^{1/p_{0}(\cdot )}%
\big\|_{p_{0}(\cdot )}$ is bounded by%
\begin{equation*}
c\Big\|\Big(\frac{g(\cdot )}{\big\|\lambda \big\|_{f_{p(\cdot ),p(\cdot
)}^{\alpha (\cdot )}}}\Big)^{p(\cdot )/p_{0}(\cdot )}\Big\|_{p_{0}(\cdot
)}+c.
\end{equation*}%
This term is bounded, since 
\begin{equation*}
\int_{\mathbb{R}^{n}}\Big(\frac{g(x)}{\big\|\lambda \big\|_{f_{p(\cdot
),p(\cdot )}^{\alpha (\cdot )}}}\Big)^{p(x)}dx\leq 1.
\end{equation*}%
\textbf{Estimation of }$\big\|\lambda ^{1}\big\|_{f_{p_{1}(\cdot
),p_{1}(\cdot )}^{\alpha _{1}(\cdot )}}$. Replacing $\alpha _{0}$, $p_{0}$
and $u$ by $\alpha _{1}$, $p_{1}$ and $v$, respectively and this leads to
the desired inequality. \ $\square $

Notice that this theorem can be generalized to the case $%
0<p_{0}^{+},p_{1}^{+}<\infty $.

\begin{theorem}
\label{calderon-prod2}Let $0<\theta <1$ and $1\leq q_{0},q_{1}<\infty $. Let 
$p_{0}\in C^{\log }$ with $1\leq p_{0}^{+}<\infty $ and $\alpha _{0},\alpha
_{1}\in C_{\mathrm{loc}}^{\log }$. We put%
\begin{equation*}
\frac{1}{p(\cdot )}:=\frac{1-\theta }{p_{0}(\cdot )},\text{ }\frac{1}{q}:=%
\frac{1-\theta }{q_{0}}+\frac{\theta }{q_{1}}\text{\ \ and\ \ }\alpha (\cdot
):=(1-\theta )\alpha _{0}(\cdot )+\theta \alpha _{1}(\cdot ).
\end{equation*}%
Then%
\begin{equation*}
\big(f_{p_{0}(\cdot ),q_{0}}^{\alpha _{0}(\cdot )}\big)^{1-\theta }\big(%
f_{\infty ,q_{1}}^{\alpha _{1}(\cdot )}\big)^{\theta }=f_{p(\cdot
),q}^{\alpha (\cdot )}
\end{equation*}%
holds in the sense of equivalent norms.
\end{theorem}

\textbf{Proof.}

\textit{Step 1.} We shall prove%
\begin{equation*}
\big(f_{p_{0}(\cdot ),q_{0}}^{\alpha _{0}(\cdot )}\big)^{1-\theta }\big(%
f_{\infty ,q_{1}}^{\alpha _{1}(\cdot )}\big)^{\theta }\hookrightarrow
f_{p(\cdot ),q}^{\alpha (\cdot )}.
\end{equation*}%
We suppose, that sequences $\lambda :=(\lambda _{j,m})_{j,m}$, $\lambda
^{i}:=(\lambda _{j,m}^{i})_{j,m},i=0,1$, are given and that%
\begin{equation*}
|\lambda _{j,m}|\leq |\lambda _{j,m}^{0}|^{1-\theta }|\lambda
_{j,m}^{1}|^{\theta }
\end{equation*}%
holds for all $j\in \mathbb{N}_{0}$ and $m\in \mathbb{Z}^{n}$. In $\mathrm{%
\eqref{notation1}}$ we replace $p(x)$ by $q$ and $\chi _{j,m}$ by $\chi
_{E_{Q_{j,m}}}$, with $E_{Q_{j,m}}\subset Q_{j,m}$ and $%
|E_{Q_{j,m}}|>|Q_{j,m}|/2$, we obtain $\mathrm{\eqref{Second-est}}$ with $%
q_{0}$ and $q_{1}$ in place of $p_{0}(x)$ and $p_{0}(x)$, respectively, and $%
\chi _{E_{Q_{j,m}}}$ in place of $\chi _{j,m}$. Estimate the second factor
by its $L^{\infty }$-norm and using Propositions \ref{prop1} we get \ $\big\|%
\lambda \big\|_{f_{p(\cdot ),q}^{\alpha (\cdot )}}\leq \big\|\lambda ^{0}%
\big\|_{f_{p_{0}(\cdot ),q_{0}}^{\alpha _{0}(\cdot )}}^{1-\theta }\big\|%
\lambda ^{1}\big\|_{f_{\infty ,q_{1}}^{\alpha _{1}(\cdot )}}^{\theta }$. 

\textit{Step 2.} We prove%
\begin{equation*}
f_{p(\cdot ),q}^{\alpha (\cdot )}\hookrightarrow \big(f_{p_{0}(\cdot
),q_{0}}^{\alpha _{0}(\cdot )}\big)^{1-\theta }\big(f_{\infty
,q_{1}}^{\alpha _{1}(\cdot )}\big)^{\theta }.
\end{equation*}%
Let the sequence $\lambda \in f_{p(\cdot ),q}^{\alpha (\cdot )}$ be given
with%
\begin{equation*}
\big\|\lambda \big\|_{f_{p(\cdot ),q}^{\alpha (\cdot )}}\neq 0.
\end{equation*}%
Let $\delta =-\frac{q}{q_{1}}$ and $\gamma (\cdot ):=\frac{p(\cdot )}{%
p_{0}(\cdot )}-\frac{q}{q_{0}}$. Observe that $\gamma $ is a constant
function. We follow ideas of the proof of Theorem 8.2 in Frazier and Jawerth 
\cite{FJ90}. Set 
\begin{equation*}
A_{\ell ,\gamma }:=\Big\{x\in \mathbb{R}^{n}:\Big(\frac{g(x)}{\big\|\lambda %
\big\|_{f_{p(\cdot ),q}^{\alpha (\cdot )}}}\Big)^{\gamma }>2^{\ell }\Big\},
\end{equation*}%
with $\ell \in \mathbb{Z}$. Obviously $A_{\ell +1,\gamma }\subset A_{\ell
,\gamma }$ for any $\ell \in \mathbb{Z}$. Now we introduce a (partial)
decomposition of $\mathbb{N}_{0}\times \mathbb{Z}^{n}$ by taking%
\begin{equation*}
C_{\ell }^{\gamma }:=\{(j,m):|Q_{j,m}\cap A_{\ell ,\gamma }|>\frac{|Q_{j,m}|%
}{2}\text{ \ and \ }|Q_{j,m}\cap A_{\ell +1,\gamma }|\leq \frac{|Q_{j,m}|}{2}%
\},\quad \ell \in \mathbb{Z}.
\end{equation*}%
The sets $C_{\ell }^{\gamma }$ are pairwise disjoint, i.e., $C_{\ell
}^{\gamma }\cap C_{v}^{\gamma }=\emptyset $ if $\ell \neq v$. Let us prove
that $\lambda _{j,m}=0$ holds for all tuples $(j,m)\notin \cup _{\ell
}C_{\ell }^{\gamma }$. Let us consider one such tuple $(j_{0},m_{0})$ and
let us choose $\ell _{0}\in \mathbb{Z}$ arbitrary. First suppose that $%
(j_{0},m_{0})\notin C_{\ell _{0}}^{\gamma }$, then either%
\begin{equation}
|Q_{j_{0},m_{0}}\cap A_{\ell _{0},\gamma }|\leq \frac{|Q_{j_{0},m_{0}}|}{2}%
\text{\quad or\quad }|Q_{j_{0},m_{0}}\cap A_{\ell _{0}+1,\gamma }|>\frac{%
|Q_{j_{0},m_{0}}|}{2}.  \label{Cond11}
\end{equation}%
Let us assume for the moment that the second condition is satisfied. By
induction on $\ell $ it follows%
\begin{equation}
|Q_{j_{0},m_{0}}\cap A_{\ell +1,\gamma }|>\frac{|Q_{j_{0},m_{0}}|}{2}\quad 
\text{for all\quad }\ell \geq \ell _{0}.  \label{Cond21}
\end{equation}%
Let $D:=\cap _{\ell \geq \ell _{0}}Q_{j_{0},m_{0}}\cap A_{\ell +1,\gamma }$.
The family $\{Q_{j_{0},m_{0}}\cap A_{\ell ,\gamma }\}_{\ell }$ is a
decreasing family of sets, i.e., $Q_{j_{0},m_{0}}\cap A_{\ell +1,\gamma
}\subset Q_{j_{0},m_{0}}\cap A_{\ell ,\gamma }$. Therefore, in view of $%
\mathrm{\eqref{Cond21}}$, the measure of the set $D$ is larger than or equal
to $\frac{|Q_{j_{0},m_{0}}|}{2}$. Hence $\varrho _{p(\cdot )}\big(\frac{g}{%
\big\|\lambda \big\|_{f_{p(\cdot ),q}^{\alpha (\cdot )}}}\big)$ is greater
than or equal to%
\begin{equation*}
\int_{\mathbb{R}^{n}}\Big(\frac{g(x)}{\big\|\lambda \big\|_{f_{p(\cdot
),q}^{\alpha (\cdot )}}}\Big)^{p(x)}\chi _{Q_{j_{0},m_{0}}\cap A_{\ell
+1,\gamma }}(x)dx\geq 2^{\ell \frac{p^{-}}{\gamma }}|D|,\quad \ell \geq \max
(\ell _{0},0).
\end{equation*}%
Now $\varrho _{p(\cdot )}\big(\frac{g}{\big\|\lambda \big\|_{f_{p(\cdot
),q}^{\alpha (\cdot )}}}\big)$ is finite, since $\lambda \in f_{p(\cdot
),q}^{\alpha (\cdot )}$, letting $\ell $ tend to infinity and using that $%
|D|\geq \frac{|Q_{j_{0},m_{0}}|}{2}$, we get a contradiction. Hence, we have
to turn in $\mathrm{\eqref{Cond11}}$ to the situation where the first
condition is satisfied. We claim,%
\begin{equation*}
|Q_{j_{0},m_{0}}\cap A_{\ell ,\gamma }|\leq \frac{|Q_{j_{0},m_{0}}|}{2}\quad 
\text{for all\ }\ell \in \mathbb{Z}.
\end{equation*}%
Obviously this yields%
\begin{equation}
|Q_{j_{0},m_{0}}\cap A_{\ell ,\gamma }^{c}|>\frac{|Q_{j_{0},m_{0}}|}{2}\quad 
\text{for all\ }\ell \in \mathbb{Z},  \label{Cond2.11}
\end{equation}%
again\ this follows by induction on $\ell $ using $(j_{0},m_{0})\notin \cup
_{\ell }C_{\ell }^{\gamma }$. Set $E=\cap _{\ell \geq \max (0,-\ell
_{0})}Q_{j_{0},m_{0}}\cap A_{-\ell ,\gamma }^{c}=\cap _{\ell \geq \max
(0,-\ell _{0})}h_{\ell }$. The family $\{h_{\ell }\}_{\ell }$ is a
decreasing family of sets, i.e., $h_{\ell +1}\subset h_{\ell }$. Therefore,
in view of $\mathrm{\eqref{Cond2.11}}$, the measure of the set $E$ is larger
than or equal to $\frac{|Q_{j_{0},m_{0}}|}{2}$. By selecting a point $x\in E$
we obtain%
\begin{equation*}
2^{j_{0}\alpha (x)}|\lambda _{j_{0},m_{0}}|\leq g(x)\leq \big\|\lambda \big\|%
_{f_{p(\cdot ),q}^{\alpha (\cdot )}}2^{-\ell /\gamma }.
\end{equation*}%
Now, for $\ell $ tending to $+\infty $ the claim, namely $\lambda
_{j_{0},m_{0}}=0$, follows. We put, as in the proof of Theorem \ref%
{calderon-prod1},%
\begin{equation*}
u(\cdot ):=q\theta \big(\tfrac{\alpha _{1}(\cdot )}{q_{0}}-\tfrac{\alpha
_{0}(\cdot )}{q_{1}}\big)+\tfrac{n}{2}\big(\tfrac{q}{q_{0}}-1\big)
\end{equation*}%
and 
\begin{equation*}
v(\cdot ):=q(1-\theta )\big(\tfrac{\alpha _{0}(\cdot )}{q_{1}}-\tfrac{\alpha
_{1}(\cdot )}{q_{0}}\big)+\tfrac{n}{2}\big(\tfrac{q}{q_{1}}-1\big).
\end{equation*}%
If $(j,m)\notin \cup _{\ell \in \mathbb{Z}}C_{\ell }^{\gamma }$, then we
define $\lambda _{j,m}^{0}=\lambda _{j,m}^{1}=0$. Let $(j,m)\in C_{\ell
}^{\gamma }$. We set 
\begin{equation*}
\lambda _{j,m}^{0}:=\lambda _{j,m,u,q_{0}}^{0,\frac{\delta }{\gamma }%
}:=2^{\ell +ju(x_{j,m})}\Big(\tfrac{|\lambda _{j,m}|}{\big\|\lambda \big\|%
_{f_{p(\cdot ),q}^{\alpha (\cdot )}}}\Big)^{q/q_{0}}.
\end{equation*}%
Also, set%
\begin{equation*}
\lambda _{j,m}^{1}:=\lambda _{j,m,v,q_{1}}^{1,\frac{\delta }{\gamma }%
}:=2^{\ell \frac{\delta }{\gamma }+jv(x_{j,m})}\Big(\tfrac{|\lambda _{j,m}|}{%
\big\|\lambda \big\|_{f_{p(\cdot ),q}^{\alpha (\cdot )}}}\Big)^{q/q_{1}}.
\end{equation*}%
Observe that%
\begin{eqnarray*}
|\lambda _{j,m}| &\leq &\big\|\lambda \big\|_{f_{p(\cdot ),q(\cdot
)}^{\alpha (\cdot )}}2^{-\ell (1-\theta +\frac{\delta }{\gamma }\theta
)}2^{-j(u(x_{j,m})(1-\theta )+v(x_{j,m})\theta )}\big(\lambda _{j,m}^{0}\big)%
^{1-\theta }\big(\lambda _{j,m}^{1}\big)^{\theta } \\
&=&\big\|\lambda \big\|_{f_{p(\cdot ),q(\cdot )}^{\alpha (\cdot )}}\big(%
\lambda _{j,m}^{0}\big)^{1-\theta }\big(\lambda _{j,m}^{1}\big)^{\theta },
\end{eqnarray*}%
which holds now for all pairs $(j,m)$. As in the proof of Theorem \ref%
{calderon-prod1}, We will prove the following two inequalities%
\begin{equation*}
\big\|\lambda ^{0}\big\|_{f_{p_{0}(\cdot ),q_{0}}^{\alpha _{0}(\cdot
)}}\lesssim 1\text{\quad and\quad }\big\|\lambda ^{1}\big\|_{f_{\infty
,q_{1}}^{\alpha _{1}(\cdot )}}\lesssim 1.
\end{equation*}%
\textbf{Estimation of }$\big\|\lambda ^{0}\big\|_{f_{p_{0}(\cdot
),q_{0}}^{\alpha _{0}(\cdot )}}$. We write%
\begin{equation*}
\sum_{j=0}^{\infty }\sum\limits_{m\in \mathbb{Z}^{n}}\Big(2^{j(\alpha
_{0}(\cdot )+\frac{n}{2})}|\lambda _{j,m}^{0}|\chi _{j,m}\Big)%
^{q_{0}}=\sum_{\ell =-\infty }^{\infty }\sum\limits_{(j,m)\in C_{\ell
}^{\gamma }}\cdot \cdot \cdot =:I.
\end{equation*}%
We use the local log-H\"{o}lder continuity of $\alpha _{0}$ and $\alpha _{1}$
to show that%
\begin{equation*}
2^{ju(x_{j,m})}\leq c\text{ }2^{ju(t)}\quad \text{and}\quad 2^{j\alpha
_{0}(x)}\leq c\text{ }2^{j\alpha _{0}(t)}\text{, \ \ }x,t\in Q_{j,m},
\end{equation*}%
where $c>0$ is independent of $\ell $ and $j$. Therefore, 
\begin{eqnarray*}
2^{j(\alpha _{0}(x)+\frac{n}{2})+ju(x_{j,m})} &\lesssim &\frac{1}{%
|Q_{j,m}\cap A_{\ell ,\gamma }|}\int_{Q_{j,m}\cap A_{\ell ,\gamma
}}2^{j(\alpha _{0}(t)+\frac{n}{2})+ju(t)}dt \\
&\lesssim &\frac{1}{|Q_{j,m}|}\int_{Q_{j,m}\cap A_{\ell ,\gamma
}}2^{j(\alpha _{0}(t)+\frac{n}{2})+ju(t)}dt \\
&\lesssim &\eta _{j,h}\ast 2^{j(\alpha _{0}(\cdot )+\frac{n}{2})+ju(\cdot
)}\chi _{Q_{j,m}\cap A_{\ell ,\gamma }}(x),\text{ \ \ }j,m\in C_{\ell
}^{\gamma },
\end{eqnarray*}%
where $h>n$ and the implicit positive constant not depending on $x$, $\ell ,m
$ and $j$. Hence $\big\|\left( I^{q_{0}}(\cdot )\right) ^{1/q_{0}}\big\|%
_{p_{0}(\cdot )}$ is bounded by%
\begin{equation*}
c\Big\|\Big(\sum_{\ell =-\infty }^{\infty }\sum\limits_{(j,m)\in C_{\ell
}^{\gamma }}\left( \eta _{j,h}\ast \digamma _{j,\ell ,m}\right) ^{q_{0}}\Big)%
^{1/q_{0}}\Big\|_{p_{0}(\cdot )},
\end{equation*}%
where 
\begin{equation*}
\digamma _{j,\ell ,m}=2^{j(\alpha _{0}(\cdot )+\frac{n}{2})+ju(\cdot )+\ell }%
\Big(\tfrac{|\lambda _{j,m}|}{\big\|\lambda \big\|_{f_{p(\cdot ),q(\cdot
)}^{\alpha (\cdot )}}}\Big)^{\frac{q}{q_{0}}}\chi _{Q_{j,m}\cap A_{\ell
,\gamma }}.
\end{equation*}%
Applying Lemma\ \ref{DHRlemma} to estimate the last norm by%
\begin{eqnarray*}
&&c\Big\|\Big(\sum_{\ell =-\infty }^{\infty }\sum\limits_{(j,m)\in C_{\ell
}^{\gamma }}\digamma _{j,\ell ,m}^{q_{0}}\Big)^{1/q_{0}}\Big\|_{p_{0}(\cdot
)} \\
&=&c\Big\|\Big(\sum_{\ell =-\infty }^{\infty }\sum\limits_{(j,m)\in C_{\ell
}^{\gamma }}2^{j(\alpha (\cdot )+\frac{n}{2})q+\ell q_{0}}\Big(\tfrac{%
|\lambda _{j,m}|}{\big\|\lambda \big\|_{f_{p(\cdot ),q}^{\alpha (\cdot )}}}%
\Big)^{q}\chi _{Q_{j,m}\cap A_{\ell ,\gamma }}\Big)^{1/q_{0}}\Big\|%
_{p_{0}(\cdot )},
\end{eqnarray*}%
since $u(\cdot )+\alpha _{0}(\cdot )=\alpha (\cdot )\frac{q}{q_{0}}+\frac{n}{%
2}\Big(\frac{q}{q_{0}}-1\Big)$. Observe that%
\begin{equation*}
2^{\ell }\leq \Big(\frac{g(x)}{\big\|\lambda \big\|_{f_{p(\cdot ),q}^{\alpha
(\cdot )}}}\Big)^{\gamma }
\end{equation*}%
for any $x\in Q_{j,m}\cap A_{\ell ,\gamma }$ and since $\gamma +\frac{q}{%
q_{0}}=\frac{p(\cdot )}{p_{0}(\cdot )}$, then $\big\|\left( I^{q_{0}}(\cdot
)\right) ^{1/q_{0}}\big\|_{p_{0}(\cdot )}$ is bounded by%
\begin{eqnarray*}
&&c\Big\|\Big(\frac{g}{\big\|\lambda \big\|_{f_{p(\cdot ),q(\cdot )}^{\alpha
(\cdot )}}}\Big)^{\gamma }\Big(\frac{g}{\big\|\lambda \big\|_{f_{p(\cdot
),q(\cdot )}^{\alpha (\cdot )}}}\Big)^{q/q_{0}}\Big\|_{p_{0}(\cdot )} \\
&=&c\Big\|\Big(\frac{g}{\big\|\lambda \big\|_{f_{p(\cdot ),q(\cdot
)}^{\alpha (\cdot )}}}\Big)^{p(\cdot )/p_{0}(\cdot )}\Big\|_{p_{0}(\cdot )}.
\end{eqnarray*}%
Obviously, the last norm is less than or equal to one. 

\textbf{Estimation of }$\big\|\lambda ^{1}\big\|_{f_{\infty ,q_{1}}^{\alpha
_{1}(\cdot )}}$. By Proposition \ref{prop1} with $E_{Q_{j,m}}^{\ell
}=Q_{j,m}\cap A_{\ell +1,\gamma }^{c}$,

\begin{equation*}
\big\|\lambda ^{1}\big\|_{f_{\infty ,q_{1}}^{\alpha _{1}(\cdot )}}\lesssim %
\Big\|\Big(\sum_{\ell =-\infty }^{\infty }\sum\limits_{(j,m)\in C_{\ell
}^{\gamma }}2^{j(\alpha _{1}\left( \cdot \right) +\frac{n}{2})q_{1}}(\lambda
_{j,m}^{1})^{q_{1}}\chi _{E_{Q_{j,m}}^{\ell }}\Big)^{1/q_{1}}\Big\|_{\infty
}.
\end{equation*}%
Observe that%
\begin{equation*}
\lambda _{j,m}^{1}\lesssim 2^{\ell \frac{\delta }{\gamma }+jv(x)}\Big(\frac{%
|\lambda _{j,m}|}{\big\|\lambda \big\|f_{p(\cdot ),q}^{\alpha (\cdot )}}\Big)%
^{\frac{q}{q_{1}}}\leq 2^{jv(x)}\Big(\frac{g(x)}{\big\|\lambda \big\|%
_{f_{p(\cdot ),q}^{\alpha (\cdot )}}}\Big)^{-\frac{q}{q_{1}}}\Big(\frac{%
|\lambda _{j,m}|}{\big\|\lambda \big\|_{f_{p(\cdot ),q}^{\alpha (\cdot )}}}%
\Big)^{\frac{q}{q_{1}}}
\end{equation*}%
for any $x\in E_{Q_{j,m}}^{\ell }$ and any $(j,m)\in C_{\ell }^{\gamma }$.
Therefore, $\big\|\lambda ^{1}\big\|_{f_{\infty ,q_{1}}^{\alpha _{1}(\cdot
)}}\lesssim 1$.\ Hence, we complete the proof. $\ \ \square $

Notice that this theorem for $\alpha _{0}=\alpha _{1}=0$ and $q_{0}=q_{1}=2$
was proved by T. Kopaliani, \cite[Theorem 3.1]{Ko09}. This theorem can be
generalized to the case $0<p_{0}^{+},p_{1}^{+},q_{0},q_{1}<\infty $.

Suppose \ that $X_{0}$ and $X_{1}$ are Banach lattices on measure space $%
\left( \mathcal{M},\mu \right) $, and let $X=X_{0}^{1-\theta }\times
X_{1}^{\theta }$ for some $0<\theta <1$. Suppose that $X$ hus the property%
\begin{equation}
f\in X,\text{ \ }\left\vert f_{n}\left( x\right) \right\vert \leq \left\vert
f\left( x\right) \right\vert \text{, }\mu \text{-a.e., \ and \ }%
\lim_{n\rightarrow \infty }f_{n}=f,\mu \text{-}a.e.\Longrightarrow
\lim_{n\rightarrow \infty }\left\Vert f_{n}\right\Vert _{X}=\left\Vert
f\right\Vert _{X}.  \label{Ca-property}
\end{equation}%
Calder\'{o}n \cite[p. 125]{Ca64} then shows that $X_{0}^{1-\theta }\times
X_{1}^{\theta }=[X_{0},X_{1}]_{\theta }$.

Now we turn to the complex interpolation of the distribution spaces\ $%
F_{p(\cdot ),p(\cdot )}^{\alpha (\cdot )}$.

\begin{theorem}
\label{Inter1}Let $0<\theta <1$. Let $p_{0},p_{1}\in C^{\log }$ with $1\leq
p_{0}^{+},p_{1}^{+}<\infty $ and $\alpha _{0},\alpha _{1}\in C_{\mathrm{loc}%
}^{\log }$. We put%
\begin{equation*}
\frac{1}{p(\cdot )}:=\frac{1-\theta }{p_{0}(\cdot )}+\frac{\theta }{%
p_{1}(\cdot )}\text{ \ and \ }\alpha (\cdot ):=(1-\theta )\alpha _{0}(\cdot
)+\theta \alpha _{1}(\cdot ).
\end{equation*}%
Then%
\begin{equation}
\lbrack f_{p_{0}(\cdot ),p_{0}(\cdot )}^{\alpha _{0}(\cdot )},f_{p_{1}(\cdot
),p_{1}(\cdot )}^{\alpha _{1}(\cdot )}]_{\theta }=f_{p(\cdot ),p(\cdot
)}^{\alpha (\cdot )}  \label{Int1}
\end{equation}%
and%
\begin{equation}
\lbrack F_{p_{0}(\cdot ),p_{0}(\cdot )}^{\alpha _{0}(\cdot )},F_{p_{1}(\cdot
),p_{1}(\cdot )}^{\alpha _{1}(\cdot )}]_{\theta }=F_{p(\cdot ),p(\cdot
)}^{\alpha (\cdot )}  \label{Int2}
\end{equation}%
holds in the sense of equivalent norms.
\end{theorem}

\textbf{Proof. }Since $f_{p_{0}(\cdot ),p_{0}(\cdot )}^{\alpha _{0}(\cdot )}$
and $f_{p_{1}(\cdot ),p_{1}(\cdot )}^{\alpha _{1}(\cdot )}$ are Banach
lattices. Then, it suffices to use Calder\'{o}n's result to $X=f_{p(\cdot
),p(\cdot )}^{\alpha (\cdot )}$\ where the property $\mathrm{%
\eqref{Ca-property}}$ follows easily from the dominated convergence theorem.
Hence Theorem \ref{calderon-prod1} yields $\mathrm{\eqref{Int1}}$. Now $%
\mathrm{\eqref{Int2}}$ follows from $\mathrm{\eqref{Int1}}$ and Theorem \ref%
{phi-tran}. $\ \ \square $

Similarly we formulate the main statement on complex interpolation of
variable Triebel-Lizorkin spaces\ $F_{p(\cdot ),q}^{\alpha (\cdot )}$.

\begin{theorem}
\label{Inter2}Let $0<\theta <1$ and $1\leq q_{0},q_{1}<\infty $. Let $%
p_{0}\in C^{\log }$ with $1\leq p_{0}^{+}<\infty $ and $\alpha _{0},\alpha
_{1}\in C_{\mathrm{loc}}^{\log }$. We put%
\begin{equation*}
\frac{1}{p(\cdot )}:=\frac{1-\theta }{p_{0}(\cdot )},\frac{1}{q}:=\frac{%
1-\theta }{q_{0}}+\frac{\theta }{q_{1}}\text{ and }\alpha (\cdot
):=(1-\theta )\alpha _{0}(\cdot )+\theta \alpha _{1}(\cdot ).
\end{equation*}%
Then%
\begin{equation*}
\lbrack f_{p_{0}(\cdot ),q_{0}}^{\alpha _{0}(\cdot )},f_{\infty
,q_{1}}^{\alpha _{1}(\cdot )}]_{\theta }=f_{p(\cdot ),q}^{\alpha (\cdot )}
\end{equation*}%
and%
\begin{equation*}
\lbrack F_{p_{0}(\cdot ),q_{0}}^{\alpha _{0}(\cdot )},F_{\infty
,q_{1}}^{\alpha _{1}(\cdot )}]_{\theta }=F_{p(\cdot ),q}^{\alpha (\cdot )}
\end{equation*}%
holds in the sense of equivalent norms.
\end{theorem}

Using a combination of the arguments used in the proof of Theorems \ref%
{Inter1} and \ref{Inter2} and the fact that $\frac{\gamma (\cdot )}{\delta
(\cdot )}$\ is a constant function with negative values, we arrive at the
following complex interpolation of variable Triebel-Lizorkin spaces.

\begin{theorem}
\label{calderon-prod11 copy(1)}Let $0<\theta <1$. Let $%
p_{0},p_{1},q_{0},q_{1}\in C^{\log }$ with $1\leq
p_{0}^{+},q_{0}^{+},p_{1}^{+},q_{1}^{+}<\infty $\ and $\alpha _{0},\alpha
_{1}\in C_{\mathrm{loc}}^{\log }$. We put%
\begin{equation*}
\frac{1}{p(\cdot )}:=\frac{1-\theta }{p_{0}(\cdot )}+\frac{\theta }{%
p_{1}(\cdot )},\text{ }\frac{1}{q(\cdot )}:=\frac{1-\theta }{q_{0}(\cdot )}+%
\frac{\theta }{q_{1}(\cdot )},\text{ }\alpha (\cdot ):=(1-\theta )\alpha
_{0}(\cdot )+\theta \alpha _{1}(\cdot )
\end{equation*}%
and%
\begin{equation*}
\gamma (\cdot ):=\frac{p(\cdot )}{p_{0}(\cdot )}-\frac{q(\cdot )}{%
q_{0}(\cdot )}.
\end{equation*}%
$\mathrm{(i)}$ We suppose that $\gamma (x)=0$ for any $x\in \mathbb{R}^{n}$.
Then%
\begin{equation*}
\lbrack f_{p_{0}(\cdot ),q_{0}(\cdot )}^{\alpha _{0}(\cdot )},f_{p_{1}(\cdot
),q_{1}(\cdot )}^{\alpha _{1}(\cdot )}]_{\theta }=f_{p(\cdot ),q(\cdot
)}^{\alpha (\cdot )}
\end{equation*}%
and%
\begin{equation*}
\lbrack F_{p_{0}(\cdot ),q_{0}(\cdot )}^{\alpha _{0}(\cdot )},F_{p_{1}(\cdot
),q_{1}(\cdot )}^{\alpha _{1}(\cdot )}]_{\theta }=F_{p(\cdot ),q(\cdot
)}^{\alpha (\cdot )}
\end{equation*}%
holds in the sense of equivalent norms.

$\mathrm{(ii)}$ We suppose that $q_{0}$ and $q_{1}$ are constants. In
addition, we assume that $\gamma (x)\neq 0$ for any $x\in \mathbb{R}^{n}$.
Then%
\begin{equation*}
\lbrack f_{p_{0}(\cdot ),q_{0}}^{\alpha _{0}(\cdot )},f_{p_{1}(\cdot
),q_{1}}^{\alpha _{1}(\cdot )}]_{\theta }=f_{p(\cdot ),q}^{\alpha (\cdot )}
\end{equation*}%
and%
\begin{equation*}
\lbrack F_{p_{0}(\cdot ),q_{0}}^{\alpha _{0}(\cdot )},F_{p_{1}(\cdot
),q_{1}}^{\alpha _{1}(\cdot )}]_{\theta }=F_{p(\cdot ),q}^{\alpha (\cdot )}
\end{equation*}%
holds in the sense of equivalent norms.
\end{theorem}

\begin{remark}
It should be mentioned that if $\alpha _{0},\alpha _{1},q_{0}$ and $q_{1}$
are constants,\ Theorem \ref{Inter-rest} is more general than what has been
given here.
\end{remark}

\textbf{Acknowledgment}

A great deal of this work has been carried out during the visit of the
author in Jena, Germany. I wish to thank Professor Winfried Sickel for the
hospitality and for many valuable discussions and suggestions\ about this
work.

\end{document}